\documentclass[12pt]{article}
\usepackage{amsmath,amsfonts,amssymb,latexsym}
\usepackage[dvips]{epsfig}
\newtheorem{theorem}{Theorem}[subsection]
\newtheorem{lemma}[theorem]{Lemma}
\newtheorem{proposition}[theorem]{Proposition}
\newtheorem{observation}[theorem]{Observation} 
\newtheorem{coro}[theorem]{Corollary}

\newcommand{\sk}{\section}
\newcommand{\ssk}{\subsection}
\newcommand{\sssk}{\subsubsection}

\newcommand{\inn}{\}_{n\in{\mathbb N}}}
\newcommand{\bull}{$\diamond$}

\newcommand{\n}{\nu}
\newcommand{\f}{\phi}
\renewcommand{\a}{\alpha}

\newcommand{\g}{\gamma}

\renewcommand{\l}{\lambda}

\newcommand{\s}{\sigma}

\newcommand{\G}{\Gamma}

\newcommand{\nb}{\bar{\nu}}
\newcommand{\begitem}{\begin{itemize}}
\newcommand{\enitem}{\end{itemize}}
\newcommand{\iti}{\item[$(i)$]}
\newcommand{\itii}{\item[$(ii)$]}
\newcommand{\itiii}{\item[$(iii)$]}
\newcommand{\itiv}{\item[$(iv)$]}
\newcommand{\itv}{\item[$(v)$]}
\newcommand{\itvi}{\item[$(vi)$]}

\newcommand{\ita}{\item[$(a)$]}
\newcommand{\itb}{\item[$(b)$]}
\newcommand{\itc}{\item[$(c)$]}
\newcommand{\itd}{\item[$(d)$]}
\newcommand{\ite}{\item[$(e)$]}
\newcommand{\itf}{\item[$(f)$]}

\newcommand{\tit}{\item[-]}

\newcommand{\preu}{{\sl Proof:~}}

\newcommand{\theo}[2]{\begin{theorem}
\label{#1}
#2
\end{theorem}}

\newcommand{\lem}[2]{
\begin{lemma}
\label{#1}
#2
\end{lemma}}

\newcommand{\cor}[2]{\begin{coro}
\label{#1}
#2
\end{coro}}

\newcommand{\obser}[2]{
\begin{observation}
\label{#1}
#2
\end{observation}}

\newcommand{\pro}[2]{
\begin{proposition}
\label{#1}
#2
\end{proposition}}

\newcommand{\nin}{n\in{\mathbb N}}

\newcommand{\m}{\mu}
\newcommand{\mb}{\bar \mu}
\newcommand{\yp}{{\mathbb H}^3}
\newcommand{\dih}{\partial_{\infty}\yp}
\newcommand{\di}{\partial_\infty M}

\newcommand{\bm}{{\bar\mu}}

\newcommand{\cpi}{{\mathbb C}{\mathbb P}^{1}}
\newcommand{\bOi}{{\cal B}^{0}_\infty}

\newcommand{\bbi}{\bar{\cal B}_{\infty}}
\newcommand{\hyp}{{\mathbb H}^{2}}

\newcommand{\auteur}{
\vskip 2truecm
\centerline{Fran\c cois Labourie}
\centerline{Topologie et Dynamique}
\centerline{Universit{\'e} Paris-Sud}
\centerline{F-91405 Orsay (Cedex)}}

\begin{document}
 \title{Random $k$-surfaces}
 \author{Fran{\c c}ois LABOURIE \thanks{L'auteur remercie l'Institut Universitaire de France.}}
 \maketitle
\sk{Introduction}
We associated in \cite{la} a compact space laminated by 2-dimensional
 leaves, to  every  compact 3-manifold $N$ with curvature less than
 -1. Considered as a ``dynamical system'', its properties generalise those of the geodesic flow. 

In this introduction, I will just sketch the construction of this space, being more precise in section \ref{introksur}. Let $k\in]0,1[$. A {\it $k$-surface} is an immersed surface in $N$, such that the product of the principal curvatures is $k$. If $N$ has constant curvature $K$, A $k$-surface has curvature $K+k$. Analytically, $k$-surfaces are described by elliptic equations. 

When dealing with ordinary differential solutions, one is lead to introduce the {\it phase space} consisting of pairs $(\gamma,x)$ where $\gamma$ is a trajectory solution of the O.D.E, and $x$ is a point on $\gamma$. We recover the dynamical picture by moving $x$ along $\gamma$.

We can mimic this construction in our situation in which a  P.D.E replaces the  O.D.E. More precisely, we can consider  the space of pairs $(\Sigma, x)$ where $\Sigma$ is a $k$-surface, and $x$ a point of $\Sigma$. 

We proved in \cite{gafa} that this construction actually makes sense. More precisely, we proved the space we just described can be compactified  by a space, called the {\it space of $k$-surfaces}. Furthermore, the boundary is finite dimensional and related in a simple way to the geodesic flow.  This space, which we denote by ${\cal N}$, is laminated by 2-dimensional leaves, in particular by those  obtained by moving $x$ along a $k$-surface $\Sigma$. A lamination means that the space has a local product structure.

The purpose of this article is to study transversal  measures, ergodic and of full support on this space of $k$-surfaces. At the present stage, let's just  to notice that since many leaves are hyperbolic ({\it cf} theorem \ref{ksur}), one cannot produce transversal  measures by Plante's argument. Our strategy will be to ``code''  by a combinatorial model on which it will be easier to build transversal  measures.

This article is organised as follows.
\begin{itemize}
\item[\bf{\ref{introksur}.}] {\bf The space of all $k$-surfaces.} We describe more precisely the {\it space of $k$-surfaces}, we are going to work with.
\item[\bf{\ref{resu}.}] {\bf Transversal  measures.} We present our main result, theorem \ref{main}, discuss other constructions and questions, and sketch the main construction.
\item[\bf{\ref{combmod}.}] {\bf A combinatorial model.} In this section,  we explain a combinatorial construction. Starting from a {\it configuration data}, we consider ``configuration spaces''.  These are  spaces of  mappings from $\mathbb{QP}^1$ to a space $W$. We produce invariant and ergodic measures under the action of $PSL(2,\mathbb Z)$ by left composition.
\item[\bf{\ref{crossratio}.}] {\bf Configuration data and the boundary at infinity of a hyperbolic 3-manifold.} We exhibit a combinatorial model associated to hyperbolic manifolds. In this context, the previous $W$ is going to be $\cpi$.
\item[\bf{\ref{coding}.}] {\bf Convex surfaces and configuration data.} We prove here that the combinatorial model constructed in the previous section actually codes for the  space of $k$-surfaces.
\item[\bf{\ref{conclu}.}] {\bf Conclusion.} We summarise our constructions and prove our main result, theorem \ref{main}.
\end{itemize}

I would like to thank W. Goldman for references about $\cpi$-structures, and R. Kenyon for discussions.
\sk{The space of all $k$-surfaces}\label{introksur}

The aim of this section is to present in a little more details the space ``of $k$-surfaces'' that we are going to work with.

Let $N$ be a compact 3-manifold with curvature less than $-1$. Let $k\in]0,1[$ be a real number. All definitions and results are expanded in \cite{la}. 

\ssk{$k$-surfaces, tubes}

If $S$ is an immersed surface in $N$, it carries several natural metrics. By definition, the {\sl $u$-metric} is the metric induced from the immersion in the unit bundle given by the Gauss map. We shall say a surface is {\it $u$-complete} if the $u$-metric is complete.

A {\it $k$-surface} is an immersed $u$-complete connected surface such that the determinant of the shape operator ({\it i.e.} the product of the principal curvatures) is constant and  equal to $k$. 

We described in \cite{la} various ways to construct $k$-surfaces. In section \ref{asym}, we summarise results of \cite{la} which allow us to obtain $k$-surfaces as solutions of an {\it asymptotic Plateau problem}. 

Since $k$-surfaces are solutions of an elliptic problem, the  germ of a $k$-surface determines the $k$-surface. It follows that a $k$-surface is determined by its image, up to coverings. More precisely, for every $k$-surface $S$ immersed by $f$ in $N$, there exists a unique $k$-surface $\Sigma$, the {\it representative of $S$}, immersed by $\phi$,  such that for every $k$-surface $\bar S$ immersed by $\bar f$ satisfying $f(S)=\bar f (\bar S)$, there exists  a covering $\pi :\bar S \rightarrow \Sigma$ such that $\bar f=\phi\circ\pi$.

By a slight abuse of language, the expression  ``$k$-surface''  will generally mean  ``representative of a $k$-surface''.

The {\it tube} of a geodesic is the set of normal vectors to this geodesic. It is a 2-dimensional submanifold of the unit bundle. 

\ssk{The space of $k$-surfaces}

The {\it space of $k$-surfaces} is  the space of pairs $(\Sigma,x)$ where $x\in \Sigma$ and $\Sigma$ is either the representative of a $k$-surface or a tube. We denote it by ${\cal N}$. It inherits a topology
coming from the topology of pointed immersed 2-manifolds in the unit bundle ({\it cf} section 2.3 of \cite{la}). Each $k$-surface (or tube)  $S_{0}$ determines a {\it leaf} ${\cal L}_{S_{0}}$ defined by
$$
{\cal L}_{S_{0}}=\{(S_{0},x)/x\in S_{0}\}.
$$

We proved in \cite{gafa} that $\cal N$ is compact. Furthermore, the partition of $\cal N$ into leaves is a lamination, {\it i.e.} admits a local product structure. Notice that ${\cal N}$ has two parts:
\begin{itemize}
\item[(1)] a  dense set which turns out to be infinite dimensional, and  which truly consists  of $k$-surfaces, 
\item[(2)] a ``boundary'' consisting of the reunion of tubes, closed,  finite dimensional, and  which is a $S^{1}$ fibre bundle over the geodesic flow. 
\end{itemize}
Therefore, in some sense, ${\cal N}$ is an extension of the geodesic flow. To enforce this analogy, one should also notice that the 1-dimensional analogue, namely  the space of curves of curvature $k$ in a hyperbolic surface is precisely the geodesic flow.

The main  theorem of \cite{la} which we quote now  shows that ${\cal N}$, as a dynamical system, enjoys the chaotic properties of the geodesic flow:

\theo{ksur}{Let $k\in ]0,1[$. Let $N$ be a compact 3-manifold. Let $h$ be a Riemannian metric on $N$  with curvature less than $-1$. Let ${\cal N}_{h}$ be the space of $k$-surface of $N$. Then
\begin{itemize}
\item[$(i)$] a generic leaf of ${\cal N}_{h}$ is dense,
\item[$(ii)$] for every positive number $g$, the union of compact leaves of ${\cal N}_{h}$ of genus greater than $g$ is dense, 
\item[$(iii)$] if ${\bar h}$ is close to $h$, then there exists a homeomorphism from ${\cal N}_{h}$ to ${\cal N}_{\bar h}$ sending leaves to leaves.
\end{itemize}}

This last property will be called the {\it stability} property.

\smallskip

To conclude this presentation, we show yet another point of view on this space, which will make it belong to a family of more familiar spaces. Assume $N$ has constant curvature, and, for just a moment,  let's vary
$k$ between $0$ and $\infty$, the range for which the associated P.D.E is elliptic. 

For $k>1$, $k$-surfaces are geodesic spheres, therefore the space of $k$-surfaces is just the unit bundle, foliated by unit spheres. 

For $k=1$, $k$-surfaces are  either  horospheres, or equidistant surfaces to a geodesic. The space of $1$-surfaces is hence described the following way:
first we take the $S^{1}$-bundle over the unit bundle, where the fibre over $u$is the set of unit vectors orthogonal $u$. This space is foliated by 2-dimensional leaves which are inverse images of geodesics. Then, we take the product of this space by $ [0,\infty[$. The number $r\in[0,\infty[$ represents the distance to the geodesic.  We complete now the space by adding horospheres, when $r$ goes to infinity. 

Our  construction  allows us to continue deforming $k$ below 1. However passing through this barrier, the space of $k$-surfaces undergoes dramatic change; in particular, it becomes infinite dimensional and ``chaotic'' as we just said.

\section{Transversal measures}\label{resu}
Let $N$ be a compact 3-manifold with curvature less than $-1$. Let $k\in]0,1[$ be a real number. Let ${\cal N}$ be the space of $k$-surfaces of $N$.

\ssk{First examples}
Let's first show  some simple examples of natural transversal  measures on ${\cal N}$. The first three ones are ergodic. They all come from the existence of natural finite dimensional subspaces in ${\cal N}$.
\begin{itemize}
\item[-] {\it Dirac measures} supported by closed leaves. By theorem \ref{ksur} $(ii)$, there are plenty of them.
\item[-] {\it Ergodic measures for the geodesic flow}. Indeed, ergodic and invariant measures for the geodesic flow give  rise to transversal  measures on the space of tubes, hence on the space of $k$-surfaces.
\item[-] {\it Haar measures for totally geodesic planes}. Assume  $N$ has constant curvature. Then, the space of oriented totally geodesic planes carries  a transverse invariant measure. Indeed, the Haar measure for $SL(2,{\mathbb C})/\pi_{1}(N)$ is invariant under the $SL(2,\mathbb R)$ action. But every oriented totally geodesic plane gives rise to a $k$-surface, namely the equidistant one to the geodesic plane. This way, we can construct an ergodic transversal  measure on ${\cal N}$, when $N$ has constant curvature. It's support is finite dimensional. 
\item[-] {\it Measures on spaces of ramified coverings}. We sketch briefly here a construction yielding transversal, but non ergodic, measures on ${\cal N}$. Let $\di$ be the boundary at infinity of the universal cover $M$ of $N$. Let $\Sigma$ be an oriented  surface of genus $g$. Let $\pi$ be a a topological ramified covering, defined  up to homeomorphism of the source,  of $\Sigma$ into $\di$. Let $S_{\pi}$ be the  set of singular points of $\pi$ and $s_{\pi}$ its cardinal. Let $S$ be a set of extra marked points of cardinal $s$. Assume $2g+s_{\pi}+s$. One can show following the ideas of the proof of theorem 7.3.3 of \cite{la} that such a ramified covering can be represented  by a $k$-surface. More precisely, there exists a solution to the asymptotic Plateau problem (as described in paragraph  \ref{asym}) represented by $(\pi,\Sigma\setminus (S_{\pi}\cup S))$. Let now $[\pi]$ be the space of ramified coverings equivalent up to homeomorphisms of the target, to $\pi$. The group $\pi_{1}(N)$ acts properly on $[\pi]$, and explicit invariant measures can be obtained using equivariant family of measures ({\it cf} section \ref{equifamili}) and configuration spaces of finite points. Since $[\pi]/\pi_{1}(N)$ is a space of leaves of $\cal{N}$, this yields transversal measures on this latter space.
\end{itemize}

None of these examples have full support, and they all have finite dimensional support.
So far, apart from these and the construction I will present in this article, I do not know of other examples of transversal  measures easy to construct.

\ssk{Main theorem} We now state our main theorem

\theo{main}{Let $N$ be a compact 3-manifold with curvature less than $-1$. Assume the metric on $N$ can be deformed, through negatively curved metrics, to a constant curvature one. Then the space of $k$-surfaces admits infinitely many mutually singular, ergodic transversal finite measures of full support.}

\ssk{First remarks}
\sssk{Restriction to the constant curvature case} The restriction upon the metric is a severe one. Actually, thanks to the stability property $(iii)$ of theorem \ref{ksur}, in order to prove our main result, it suffices to show the existence of transversal  ergodic finite measures of full support in the case of constant curvature manifolds.

\sssk{Choices made in the construction} The measure we construct  on $\cal N$ depends on several choices, and various choices lead to mutually singular measures.

We describe now one of the crucial choice  needed in the construction. 

Let $M$ be the universal cover of $N$. Let $\di$ be its boundary at infinity. Let $\cal P (\di)$ be the space of probability radon measures on $\di$. Let 
$$
{O}_{3}=\{(x,y,z)\in\di^{3}/x\not=y\not=z\not=x\}.
$$
The construction requires a map $\nu$, invariant under the natural action of $\pi_{1}(N)$,
$$
O_{3}\stackrel{\nu}{\longrightarrow} {\cal P} (\di).
property$$
Here, $\nu(x,y,z)$ is assumed to be of full support, and to fall in the same measure class, independently of $(x,y,z)$. Such maps are easily obtained through {\it equivariant family of measures} (also described in F.~Ledrappier's article \cite{led} as {\it Gibbs current}, {\it crossratios} {\it etc}) and a barycentric construction as shown in paragraph \ref{mesconfdata}. 

\ssk{Strategy of proof} As we said in the introduction, the construction is obtained through a coding of the space of $k$-surfaces. We give now a heuristic, non rigorous, outline of the proof, which is completed in the last section.

From the stability property, we can assume $N$ has constant curvature. Our first step (section \ref{coding}) is to associate to (almost) every $k$-surface a locally convex pleated surface, analogous to a ``convex core boundary''. It turns out that this way we can describe a dense subset of $k$-surfaces, by locally convex pleated surfaces, and in particular by their {\it pleating loci} at infinity. Such pleating loci are described as special maps from $\mathbb{QP}^{1}$ to $\mathbb{CP}^1$. This is the aim  of  sections \ref{crossratio} and \ref{coding}. Identifying $\mathbb{QP}^{1}$ with the space of connected components of $\hyp$ minus a trivalent tree, we build invariant measures on this space of maps as projective limit of measures on finite configuration spaces of points on $\cpi$. This is done in section \ref{combmod}.   

\ssk{Comments and questions}

\sssk{General negatively curved 3-manifolds} As we have seen before, the proof only works in the case of constant curvature manifolds, extending to other cases through the stability . Of course, it would be more pleasant to obtain transversal measures without any restriction on the metric. Some parts of the construction do not require any hypothesis on the metric, and we have tried to keep, at the price of slightly longer proof sometimes, the proof as general as possible. 
\sssk{Equirepartition of closed leaves} Keeping in mind the analogy with the geodesic flow and the construction of the Bowen-Margulis measure, here is a completely different  attempt to exhibit transversal measures, without any initial assumption on the metric. Define the {\it H-area} of a $k$-surface to be the integral of its mean curvature. It is not difficult to show that for any real number $A$, the number $N(A)$ of $k$-surfaces in ${\cal N}$ of $H$-area less than $A$ is bounded. Starting from this fact, one would like to know if closed leaves are equidistributed in some sense, {\it i.e.} that some average $\mu_{n}$ of measures supported on closed leaves of area less than $n$  weakly converges as $n$ goes to infinity. We can be more specific and ask about closed leaves of a given genus, or closed leaves whose $\pi_{1}$ surjects onto a given group. This is a whole range of questions on which I am afraid to say I have no hint of answer. However, the constructions in this article should be related to equirepartition of  ramified coverings of the boundary at infinity by spheres.

\section{A combinatorial model}\label{combmod}
In general, ${\cal P}(X)$ will denote the space of probability radon measures
on the topological space $X$,  $\delta_{x}\in {\cal P}(X)$ will
be  the Dirac measure concentrated at $x\in X$, and ${\mathbb I}_S$ will
be the characteristic function of the set $S$.

In this section, we shall describe {\it restricted infinite configuration spaces} (\ref{configspa}), which are roughly speaking spaces of infinite sets of points on a topological space $W$, associated to a {\it configuration data} (\ref{data}). Our main result is theorem \ref{exismes} which defines invariant ergodic measures of full support on these spaces, starting from measures defined on configuration data as in paragraph \ref{defconfdata}. One may think of these restricted infinite configuration spaces  as analogue of subshifts of finite type, where the analogue of the Bernouilli shift is the space of maps of $\mathbb{QP}^{1}$ (instead of $\mathbb Z$) into a space $W$ with the induced action of $PSL(2,{\mathbb Z})$. We call this latter space infinite configuration space and describe in the first paragraph, as well as related notions. The role of the {\it configuration data} is  that of local transition rules. 

\subsubsection{Trivalent tree}\label{3tree} 
\begin{figure}[h!]  
  \centerline{\psfig{figure=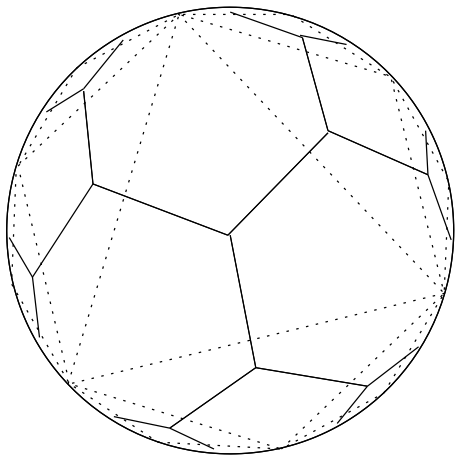}}
  \caption{\it The infinite trivalent tree dual to the ideal triangulation}\label{tree}
\end{figure}

We consider the infinite trivalent tree $T$, with a fixed
cyclic ordering on the set of edges stemming from any  vertex. Alternatively we can think of  this
ordering as defining a proper embedding of the tree in the real plane ${\mathbb R}^{2}$, such that
the cyclic ordering agrees with the orientation. Another  useful picture to
keep in mind is to consider the periodic tiling of the hyperbolic plane $\hyp$  by ideal
triangle, and our tree is the dual to this picture (Figure \ref{tree}). The group $F$ of
symmetries of that picture, which we abusively call the  {\it ideal triangle
group}, acts transitively on the set of vertices. It is isomorphic to $F={\mathbb Z}_{2}*{\mathbb Z}_{3}=PSL(2,{\mathbb Z})$.

We now consider the set $B$ of
connected components of $\hyp\setminus T$. In our tiling picture
this set $B$ is in one-to one correspondence with  the set of vertices
of triangles, and it follows that the ideal triangle group 
$F$ acts also transitively on $B$. Actually $B$ can be identified with $\mathbb{QP}^{1}$ and this identification agrees with the action of $PSL(2,{\mathbb Z})$.

\subsubsection{Quadribones, tribones}

\begin{figure}[h!]  
  \centerline{\psfig{figure=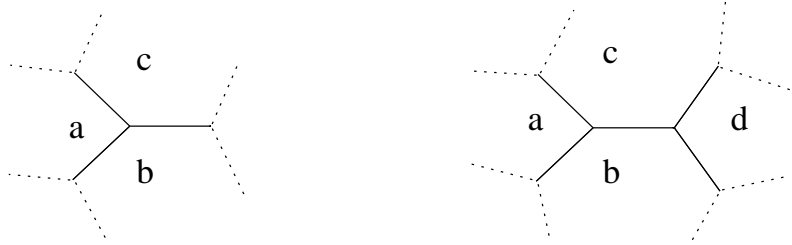}}
  \caption{\it tribone $(a,b,c)$ and quadribone $(a,b,c,d)$}\label{bone}
\end{figure}

Every edge of $T$ defines a set of four points in $B$, namely the
connected components of $\hyp\setminus T$ that touches this
edge, and we shall call these particular sets {\it quadribones}.  We
consider this set as an oriented set, {\it i.e.} up to signature 1
permutations, as labelled in the figure \ref{bone}. 
Also,
every vertex of the tree defines special subsets of three points in $B$,
that we shall call {\it tribones}. Obviously every quadribone contains
two tribones corresponding to the extremities of the corresponding edge,
and again these quadribones are oriented sets. When our quadribone is
given by $(a,b,c,d)$ the two corresponding tribones are $(a,b,c)$ and
$(d,c,b)$.

\subsubsection{Infinite configuration spaces}\label{configspa}

We define the {\it infinite configuration
space} of $W$ to be the  space, denoted  ${\cal B}_{\infty}$, of maps
from $B$ to $W$.

Notice that every tribone $t$ (resp. quadribone $q$) of $B$ defines a natural map 
from ${\cal B}_{\infty}$ to $W^{3}$ (resp. $W^{4}$) given respectively 
by  $f\mapsto f(t)$ and $f\mapsto f(q)$; we call these maps {\it associated
maps to the tribone $t$} (resp. {\it to the quadribone $q$}).

\subsection{Local rules}\label{data}
For the construction of our combinatorial model, we need the following
definitions. 
\subsubsection{Configuration data}
We shall say that $(W,\G, O_{3},O_{4})$
defines a {\it (3,4)-configuration data} if:
\begitem
\ita    $W$ is a metrisable topological space;
\itb  $\G$ is a discrete group acting continuously on $W$.
\enitem
\medskip
We deduce from that a (diagonal) action of $\G$ on $W^n$ which commutes
with the action of the $n^{th}$-symmetric group $\s_{n}$. Let $\s_{n}^{+}$ be the subgroup of $\s_{n}$ of signature $+1$. Let  $\l_{3}=\s_{3}^{+}$,  and $\l_{4}\subset\s_{4}^{+}$ be  the subgroup generated by $(a,b,c,d)\mapsto (d,c,b,d)$. Let $\Delta_{n}=\{(x_{1},\ldots,x_{n})/\exists i\not= j,~x_{i}=x_{j}\}$.
We assume
furthermore
that: 
\begitem
\itc $O_{n}$ are open $\l_{n}\times \G$-invariant subsets  of 
$W^{n}\setminus\Delta_{n}$,
on which $\G$ acts properly.

\itd $p( O_{4})=O_{3}$, where $p$ is the projection from $W^4$ to
$W^3$ defined by 
$$
(a,b,c,d))\mapsto (a,b,c).
$$

\enitem

We shall also say a configuration data is {\it Markov} if it satisfies the following extra hypothesis:

\begin{itemize}
\ite  There exist some constant $p\in {\mathbb N}$, such that if $(a,b,c)$ and
$(d,e,f)$ both belong to $O_{3}$, then there exists a sequence
 $(q_{1},\ldots,q_{j})$ of elements  of $O_{4}$, where $j\leq p$ and
 $q_{n}=(q^{1}_{n},q^{2}_{n},q^{3}_{n},q^{4}_{n})$, satisfying:
\begin{itemize}
\item $(q_{1}^{1},q_{1}^{2},q_{1}^{3})=(a,b,c)$ ;
\item $(q_{j}^{2},q_{j}^{3},q_{j}^{4})=(d,e,f)$ ;
\item $(q^{2}_{n},q^{4}_{n},q^{3}_{n})=(q^{1}_{n+1},q^{2}_{n+1},q^{3}_{n+1})$
 or $(q^{3}_{n},q^{2}_{n},q^{4}_{n})=(q^{1}_{n+1},q^{2}_{n+1},q^{3}_{n+1})$.
\end{itemize}

\enitem

In paragraph \ref{conseqe}, this property will have a geometric consequence.

\subsubsection{Measured configuration data}\label{defconfdata}
Our next goal is to associate measures to this situation. We shall say  $(W,\G, O_{3},O_{4},\m^{3},\m^{4})$ is a  {\it (3,4)-
measured configuration data} if: 
\begitem
\itf  $\mu^{n}$ are $\l_{n}\times\G$-invariant  measures, such that
$p_{*}\mu^{4}=\mu^{3}$.

\item[(g)] The pushforward measures 
on $O_{n}/\G$ are probability measures.
\enitem
\medskip

We shall say that the measured configuration data is {\it regular} if it satisfies the following extra condition:
\begitem
\item[(h)] the measure $\m_{4}$ is in the measure class of 
${\mathbb I}_{O_{4}}m\otimes m\otimes m \otimes m$ where $m$ is of full
support in $W$. It follows that  $\m_{3}$ is in the measure class of 
${\mathbb I}_{O_{3}}m\otimes m\otimes m $. 
\enitem

We also say two regular measured configuration data $(W,\Gamma,O_{3},O_{4},\m^{3},\m^{4})$ and
$(W,\Gamma,O_{3},O_{4},\bar\m^{3},\bar\m^{4})$, defined on the same configuration data, are {\it mutually singular} if $\m^{3}$ and $\bar\m^{3}$ are mutually singular.  
\subsubsection{Remarks}{\label{remdata}}

$(i)$~~ From disintegration of measures, it follows from the hypothesis (f) and (g)  that for $\mu^{3}$-almost every triple of
points $(a,b,c)$ in $W$, we have a probability measure $\nu_{(a,b,c)}$ on $W$
such that for every positive measurable function $f$ on $W^{4}$:
$$
\int_{W^{4}}f(a,b,c,d)d\mu^{4}(a,b,c,d)=\int_{W^{3}}\bigg(\int_{W}f(a,b,c,d)d\n_{(a,b,c)}(d)\bigg)d\mu^{3}(a,b,c)
$$
which we can rewrite as
$$
d\mu^{4}=\int_{W^{3}}\big(\delta_{(a,b,c)}\otimes d\n_{(a,b,c)}\big) d\mu^{3}(a,b,c),
$$
$(ii)$~~\label{datacheck} Conversely, here is a way to build a regular measured
configuration data starting from a configuration data $(W,\G, O_{3},O_{4})$, 
{\sl if we assume that $O_{4}$ is invariant under $\s_{4}^{+}$.}

Assume we have:
\begitem
\tit a $\G$-invariant measure   $\mb^{3}$ on $W^{3}$ in the measure
class of ${\mathbb I}_{O_{3}}m\otimes m\otimes m$ where $m$ has full
support, such that the push forward on $O_{3}/\Gamma$ is a probability
measure;

\tit a $\G$-equivariant map $\nb$: 
\begin{displaymath}
\left\{ \begin{array}{cc}
O_{3} & \rightarrow {\cal P}_{m}(W) \\
(a,b,c) & \mapsto\nb_{(a,b,c)} 
\end{array}
\right.
\end{displaymath}
where ${\cal P}_{m}(W)$ is the set of finite radon measures on
$W$ in the measure class of $m$.
\enitem
Then, we can build $\m^{3}$ and $\m^{4}$ which will fulfill the
requirements of the definition. Let's describe the procedure:

Firstly,  we define a probability measure $\mb^{4}$
on $O_{4}$ to be proportional to
$$
{\mathbb I}_{O_{4}}\int_{W^{3}}\big( \delta_{(a,b,c)} \otimes \nb_{(a,b,c)}\big) d\mb^{3}(a,b,c).
$$

Secondly,
we average $\mb^{4}$ using the group $\s_{4}^{+}$  and  obtain a finite
measure $\m^{4}$ on $O_{4}/\G$, and we ultimately take $\m^{3}=p_{*}\m^{4}$.

It is a routine check now that $\m^{3}$ and $\m^{4}$ defined this way
satisfy our needs.

Furthermore, if $\mb^{3}$ has full support in $O_{3}$ as well as
$\nu_{(a,b,c)}$  for
$\mb^{3}$-almost every $(a,b,c)$ in $W_{3}$, then $\mu^{3}$ and
$\mu^{4}$ have full support.

\subsubsection{Example}

\label{mainexample} In the sequel, we only wish to study one example that we describe briefly now and more precisely in section \ref{crossratio}. Our specific interest lies in the following situation.
\begitem
\tit  $\G$ is a cocompact discrete subgroup of $PSL(2, {\mathbb C})$;
\tit  $W=\cpi$ with the canonical action
of $\G$; it is a well known fact that  $\G$ acts properly
on 
$$
U_{n}=\{(x_{1},\ldots,x_{n})\in( \cpi
 )^{n}/~x_{i}\not=x_{j}~{\rm if}~i\not=j\}.
$$
\tit $O_{3}=U_{3}$,  
\tit $O_{4}$ is  the set of points whose crossratios have a non zero
imaginary part; it will  satisfy hypothesis $(e)$ for $N=1000$ ({\it cf} paragraph \ref{crossratio2}).
\enitem

This is a Markov configuration data and furthermore in this specific situation $O_{4}$ is is invariant under $\s_{4}^{+}$. We will explain in paragraph \ref{mesconfdata}  how to attach measures to this situation, and discuss the case of general negatively curved 3-manifolds.

\subsubsection{Final remark}
Even though we only wish to study that specific class of examples, it is
a little more comfortable to work  in our more general setting, since very little
of the geometry is used at this stage.

\subsection{Restricted infinite configuration spaces and main result}

Let now $(W,\G,O_{i})$ be a (3,4)-configuration data  ({\it cf} \ref{data}),

We define the {\it restricted infinite  configuration
space} of $W$ to be the subset  $\bar{\cal B}_{\infty}$ of
${\cal B}_{\infty}$, consisting of those maps such that the image of every tribone lies in
$O_{3}$, and the image of every quadribone is in $O_{4}$.
$$
\bbi=\{f\in {\cal B}_{\infty}/ {\rm for~all~tribone}~t, {\rm quadribone}~q, f(t)\in O_{3}, f(q)\in O_{4}\}.
$$

Let also ${\cal B}^0 _{\infty}$ be the open set of the infinite
configuration space such the image of at least one tribone lies in
$O_{3}$. Let's call this subset the {\it non degenerate configuration
space}, and notice that $\Gamma$ acts properly on this open subset of
${\cal B}_\infty$.

Now we can state the theorem we wish to prove:

\theo{exismes}{Let  $(W,\G,O_{i},\m^{i})$ be a (3,4)-measured configuration
data. Then there exists a $\G$-invariant measure $\mu$ on the infinite 
configuration space of $W$, which is invariant by the action of the ideal 
triangle group, such that: 
\begin{itemize}
\iti the restricted infinite configuration
space $\bar{\cal B}_\infty$ is of full measure and $\mu$ has full support on it provided the data is regular;
\itii the pushforward of $\mu$ on  ${\cal B}^0 _{\infty}/\G$ is finite,
where ${\cal B}^0 _{\infty}$ is the non degenerate infinite configuration
space;
\itiii given any tribone or quadribone, the pushforward of
$\mu$ by the associated maps on $W^{3}$ and $W^{4}$ are our original
 $\mu^{3}$, $\mu^{4}$;
\itiv two regular mutually singular measured configuration data give rise to
mutually singular measures ;
\itv if the configuration data is Markov and regular, then the pushforward of $\mu$  on  ${\cal B}^0 _{\infty}/\G$  is ergodic with respect to the action of the
infinite triangle group.
\end{itemize} 
}

Essentially, this measure is built by a Markov type procedure.

\subsection{Construction of the measure} Let  $(W,\G,O_{i},\m^{i})$ be a (3,4)-configuration
data. We shall use the notations and definitions of the preceding
sections. 

Also in our constructions, for every $(x,y,z)\in O_{3}$, we shall denote by  $\n_{(x,y,z)}$  the probability
measure coming from the disintegration of $\m^{4}$ over $\m^{3}$ as defined in
paragraph \ref{remdata}.

\subsubsection{Connected sets, $P$-bones, $P$-disconnected sets}\label{notatree}
For our constructions, we need  a terminology for some subsets of $B$ which roughly correspond to certain subtrees of $T$.

\smallskip
A subset $A$ of $B$ will be called {\it connected} if it is a reunion of
quadribones such that the reunion $e(A)$ of the associated edges is
connected; if $v$ is a vertex, it will be called {\it $v$-connected} if
furthermore $e(A)$ contains $v$. In other words a connected subset of
$B$ is the reunion of the connected components touching the edges of a
connected subtree of $T$.

\smallskip
A subset $A$ of $B$ will be called a {\it $P$-bone} if it is connected and
the reunion of less than $P$ quadribones; two subsets $A$ and $C$ will
be called {\it $P$-disconnected} if there is no $P$-bone which
intersects both $A$ and $C$. 

\subsubsection{Relative configuration spaces}
If $A$ is a subset of $B$, we shall note: 
\begin{itemize}
\tit ${\cal W}(A)$ the set of maps from $A$ to $W$; in particular,
${\cal W}(B)={\cal B}_{\infty}$.

\tit $\bar{\cal W}(A)$ the set of maps such that the image of every tribone
of $A$ 
lies in $O_{3}$, and the image of every quadribone is in $O_{4}$; if $A$ is finite, $\bar{\cal W}(A)$ is an open 
set on which $\G$ acts properly. Again,
$\bar{\cal W}(B)=\bar{\cal B}_{\infty}$.
\end{itemize}

\subsubsection{Finite construction}
We can now prove:  
\pro{finitemes}{Let $A$ be a  finite $v_{0}$-connected subset of $B$. Then,
there exists a radon measure $\mu^{A,v_{0}}$ on  $\bar{\cal W}(A)$
enjoying the following properties:
\begin{itemize}
\iti the pushforward of $\mu^{A,v_{0}}$ on $\bar{\cal W}(A)/\G$ is
finite ; it is  of full support if the data is regular;
\itii let $t_{0}$ be the  tribone corresponding to the vertex $v_{0}$;
let also $t_{0}$ be the associated  map from
$A$ to $W^{3}$; then $t_{*}\mu^{A,v_{0}}=\mu^{3}$;
\itiii let $q$ be  a $v_{0}$-connected quadribone ; assume $q\subset A$;  let $q$ be the associated map from $A$ to $W^{4}$; then  $q_{*}\mu^{A,v_{0}}=\mu^{4}$
\itiv assume there exist a tribone $t\subset A$, some element $a\in
B\setminus A$, such that $q=t\cup\{a\}$ is a quadribone ; let now $C=A\cup \{a\}$
and identify ${\cal W}(C)$ with ${\cal W}(A)\times W$ then
$$
\mu^{C,v_{0}}=\int_{{\cal W}(A)}\big(\delta_{f}\otimes \nu_{f(t)}\big)
d\mu^{A,v_{0}}(f).
$$ 
\itv let  $A\subset C$ ; let $p$ be  the natural restriction from $\bar{\cal W}(C)$ to $\bar{\cal W}(A)$. Then $p_{*}\mu^{C,v_{0}}=\mu^{A,v_{0}}$ ;
\itvi if $(\mu^{3},\mu^{4})$ and $(\bm^{3},\bm^{4})$ are regular and mutually singular, then the corresponding measures $\mu^{A,v_{0}}$ and $\bm^{A,v_{0}}$ are mutually singular.
\end{itemize}}

One should notice that the listed properties defines $\mu^{A,v_{0}}$
uniquely. 

We shall also say in the sequel that if $C$ and $A$ are as in $(iv)$,
that $C$ is obtained from $A$ by {\it gluing a quadribone along a tribone}, as in figure \ref{glu}.

\begin{figure}[h!]
 \centerline{\psfig{figure=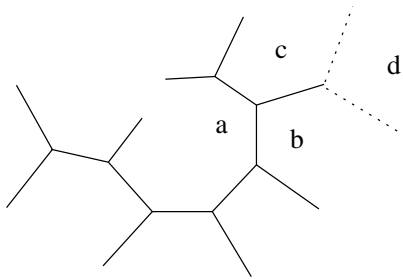}}
 \caption{\it Gluing a quadribone $(a,b,c,d)$ along a tribone $(a,b,c)$ }\label{glu}
\end{figure}

We have a useful consequence of the previous proposition:
\cor{fininv}{Let $A$ be a finite set and let $v$ and $w$ such that $A$
is both $v$-connected and $w$-connected, then $\mu^{A,v}=\mu^{A,w}$.}

Now of course, we may write $\mu^{A}=\mu^{A,v}$.

Our last proposition exhibits some kind of 
``Markovian'' property of our measure. 
\pro{markov}{Assume the configuration data is Markov and regular. There exists an integer $P$, such that if  $A_{0}$ and
$A_{1}$ are two $P$-disconnected subsets of a finite set $C\subset B$, then  $
(p^{0},p^{1})_{*}\mu^{C}$ and $p_{*}^{0}\mu^{C}\otimes p_{*}^{1}\mu^C$ are in the same measure class. Here, 
$p^{i}~:
{\cal W}(C)\rightarrow {\cal W}(A_{i})$ are the natural restriction maps. }

We will now prove the results stated in this paragraph.

\subsubsection{Proof of proposition \ref{finitemes}}\label{prooffinitemes}
Introduce first some notations with respect to a vertex $v$. By definition $B_{n}(v)$ will denote the union of all $v$-connected $n$-bones; also, for any subset $A$ of $B$,  we pose $A_{n}(v)=B_{n}(v)\cap A$.

For the moment, we will work with a fixed $v_{0}$ and will omit the dependence in
$v_{0}$ in the notations for the sake of simplicity, in particular $A_{n}=A_{n}(v_{0})$. We will construct this measure by an induction procedure.

Our first task is to build  for every $n\in{\mathbb N}$, a  map:
\begin{displaymath}
\n^{A,n}:\left\{
\begin{array}{rcl}
  {\bar{\cal W}}(A_{n}) & \rightarrow  & {\cal P}\big({\cal W}(A_{n+1}\setminus
A_{n})\big)\\
   f &\mapsto &\nu_{f}^{A,n}.\\
\end{array}
\right.
\end{displaymath}
Let's do it. If $a\in A_{n+1}\setminus A_{n}$, it belongs to a unique quadribone
$q_a \subset A_{n+1}$. Let 
$t_{a}=q_{a}\setminus\{a\}$; notice that $t_{a}$ is a subset of $A_{n}$.
Let $A_{n+1}\setminus A_{n}=\{a_{1},\ldots,a_{q}\}$. In particular,
${\cal W}(A_{n+1}\setminus A_{n})$ is identified with $W^{q}$.
Let $T_{n}^{A}=\cup_{i=1}^{i=q}t_{a_{i}}$. We have a natural restriction map
$$
i^{A,n}:~{\bar{\cal W}}(A_{n})\longrightarrow{\bar{\cal W}}(T_{n}^{A}).
$$
We define
\begin{displaymath}
{\bar\nu}^{A,n}:\left\{
\begin{array}{rcl}
  {\bar{\cal W}}(T_{n}^{A}) & \rightarrow  &{\cal P}( W^{q})={\cal P}\big({\cal W}(A_{n+1}\setminus
A_{n})\big)\\
   f &\mapsto &\bigotimes_{i}\nu_{f(t_{i})}.\\
\end{array}
\right.
\end{displaymath}
Finally, we set: $\n^{A,n}={\bar\nu}^{A,n}\circ i^{A,n}$.

Next, notice the following fact.
Let $f\in {\bar{\cal W}}(A_{n})$. Let ${\bar{\cal W}}_{f}(A_{n+1})$ be the fibre, over $f$, of the restriction map. Let's use  the identification
$$
{\cal W}(A_{n+1})={\cal W}(A_{n+1}\setminus A_{n})\times{\cal
W}(A_{n}).
$$
Then, ${\bar{\cal W}}_{f}(A_{n+1})$ has full measure  for
$\n^{A,n}_{f}\otimes \delta_{f}$.

We can now define our measure on ${\bar{\cal W}}(A_{n+1})$ by an induction procedure:
\begin{itemize}
\tit  ${\bar{\cal W}}(A_{0})$ is identified with $O_{3}$ using $t_{0}$;  we
define $\m^{A_{0}}=(t_{0}^{-1})_{*}\m^{3}$;
\tit  assuming by induction that  $\m^{A_{n}}$ is defined on 
${\cal W}(A_{n})$ such that ${\bar{\cal W}}(A_{n})$ has full measure, we set
$$
\m^{A,{n+1}}=\int_{{\bar{\cal W}}(A_{n})}\big(\n^{A,n}_{f}\otimes
\delta_{f}\big) d\mu^{A,n}(f).
$$
\end{itemize}
From the previous observation, we deduce that  ${\bar{\cal W}}(A_{n+1})$ has full measure. Furthermore, if the $\mu^{i}$ have full support, then
$\m^{A,{n+1}}$ has full support.

Finally, there exists $p\in{\mathbb N}$ such that $A=A_{p}$, and we define
$$
\mu^{A,v_{0}}=\mu^{A,p}.
$$
Properties $(i)$, $(ii)$, $(iii)$, and $(vi)$ are immediately checked.
Let's prove property $(iv)$. 

Notice first that  $a$ lies in exactly one quadribone $q$ of $C$. Let 
$d$ be the unique tribone of $C$ that contains $a$. 
Then, there exists  $p_{0}$ such that
\begin{eqnarray*}
C_{p}&=&A_{p}\hbox{~for~}p<p_{0},\\
C_{p}&=&A_{p}\cup\{a\}\hbox{~for~}p\geq p_{0}.
\end{eqnarray*}
By construction, using the obvious identifications, we have
\begin{eqnarray*}
\mu^{C,p}&=&\mu^{A,p},\hbox{~~for~}p<p_{0},\\
\mu^{C,p}&=&\int_{{\cal W}(A_{p})}(\delta_{f}\otimes\nu_{f(q\setminus a)})d\mu^{A,p}(f),\hbox{~~for~}p=p_{0}.~~~(*)
\end{eqnarray*}
To conclude the proof of $(iv)$, it remains to prove $(*)$ for $p>p_{0}$. By induction, this follows from the fact that, for $p>p_{0}$, $T^{A}_{p}=T^{C}_{p}$.Let's check that step by step.
By definition,  
$$
\m^{C,{p+1}}=\int_{{\cal W}(C_{p})}\big(\n_{f(T_{p})}\otimes
\delta_{f}\big) d\mu^{C,p}(f).
$$
But, by induction
$$
\mu^{C,p}=\int_{{\cal W}(A_{p})}(\delta_{g}\otimes\nu_{g(q\setminus a)})d\mu^{A,p}(g).
$$
Combining the two last equalities, and using $T^{A}_{p}=T^{C}_{p}$, we get
\begin{eqnarray*}
\m^{C,{p+1}}&=&\int_{{\cal W}(A_{p})}(\delta_{g}\otimes\nu_{g(q\setminus a)}\otimes \nu_{g(T_{p}^{A})})d\mu^{A,p}(g)\\
&=& \int_{{\cal W}(A_{p+1})}(\delta_{g}\otimes\nu_{g(q\setminus a)})d\mu^{A,p+1}(g) .
\end{eqnarray*}
This is what we wanted to prove.

Property $(v)$ is an immediate consequence of $(iv)$. Indeed, if $C$ contains
$A$, it is obtained inductively from $A$ 
by ``gluing quadribones along tribones'' as in $(v)$. \bull
\subsubsection{Proof of corollary \ref{fininv} } 

Obviously, it suffices to
prove it whenever $v$ and $w$ are the extremities of a common edge $e$. Let $q$ be the associated quadribone. Since we can build $A$
from $q$ by successively ``gluing quadribones along tribones'', using property $(v)$ of
proposition \ref{finitemes}, it suffices to show that
$$
\mu^{q,v}=\mu^{q,w}.
$$ 
Thanks to proposition \ref{finitemes} $(iii)$, this follows from the invariance of
$\mu^{4}$
under the permutation $(a,b,c,d)\mapsto (d,c,b,a)$.\bull

\subsubsection{A consequence of hypothesis (e) of \ref{data}}
Using the previous notations, we have:
\pro{conseqe}{Assume the configuration data is Markov. Then, there exists an integer $P$,
such that if $A_{0}$ and $A_{1}$ are connected and $P$-disconnected, if $C$ is a connected set that contains both then
$$
(p^{0},p^{1})(\bar{\cal W}(C))=\bar{\cal W}(A_{0})\times
\bar{\cal W}(A_{1}).
$$}
\preu Let  $A_{0}$ and $A_{1}$ be two $P$-disconnected subsets. Then there exists a $N$-bone $K$, where $N>P$, such that $K$ intersects each $A_{i}$ exactly  along one tribone $t_{i}$ as in figure \ref{conn}. Let $D=A_{0}\cup K\cup A_{1}$.

\begin{figure}[h!]
\centerline{\psfig{figure=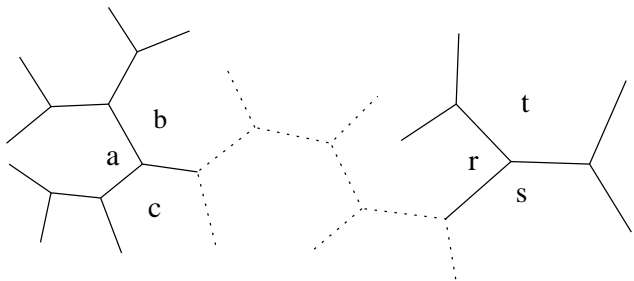}}
\caption{\it $t_{0}=(a,b,c)$, $t_{1}=(r,s,t)$ }\label{conn}
\end{figure}

Let $f_{0}$ (resp. $f_{1}$) be an element of $\bar{\cal W}(A_{0})$ (resp. $\bar{\cal W}(A_{1})$). Hypothesis (e) of \ref{data} implies there exists some element $g$ of $\bar{\cal W}(K)$ such that $g$ coincides with $f_{0}$ (resp. $f_{1}$) on $t_{0}$ (resp $t_{1}$) provided $N$ is greater than $p$. Gluing together $g$ and the $f_{i}$, we obtain an element $h$ of $\bar{\cal W}(D)$, whose restriction to $A_{i}$ is $f_{i}$. In other words, the restriction from ${\bar{\cal W}}(D)$ to ${\bar{\cal W}}(A_{0})\times{\bar{\cal W}}(A_{1})$ is surjective. To conclude,  its suffices to notice that since $D$ is connected, the restriction from $\bar{\cal W}(C)$ to  $\bar{\cal W}(D)$ is surjective.\bull

\subsubsection{Proof of proposition  \ref{markov}}
The first  point to notice is that if $\m^{3}$ and $\m^{4}$ are in the
measure class of ${\mathbb I}_{O_{3}}m\otimes m\otimes m$ and ${\mathbb I}_{O_{4}}m\otimes m\otimes m\otimes
m$ respectively, then for $m$-almost every tribone $t$, $\n_{t}$ is also
in the measure class of ${\mathbb I}_{O_{t}}m$ where $O_{t}$ is such that
$\{ t \}\times O_{t} =p^{-1}(t)\cap O_{4}$. 
\smallskip
It follows that 
if $A$ is connected then $\m^{A}$ is actually in the measure class of
$$
{\mathbb I}_{\bar{\cal W} (A)}m^{\otimes \# A}.
$$ 
Let's prove this last assertion by induction: let $C_{n}$, $1\leq n\leq p$  be an increasing
sequence of sets such that $C_{1}$ is  quadribone, $C_{p}=A$ and
$C_{n+1}$ is obtained from $C_{n}$ by gluing a quadribone along a
tribone $t_{n}$ as in proposition \ref{finitemes} $(iv)$. Then 
$$
\bar{\cal W}(C_{n+1})=\bigcup_{f\in \bar{\cal W}(C_{n})}\{f\}\times O_{f(t_{n})},
$$
where we identified ${\cal W}(C_{n+1})$ et ${\cal W}(C_{n})\times
W$. An inductive use of  \ref{finitemes} $(iv)$ implies our statement. 

\smallskip
Assume now the configuration data is Markov. Then according to proposition 
\ref{conseqe}, we have
$$
(p^{0},p^{1})(\bar{\cal W}(C))=\bar{\cal W}(A_{0})\times
\bar{\cal W}(A_{1}).
$$
Hence our proposition.\bull

\subsubsection{Infinite construction, and proof of properties  $ (i),\ldots ,(iv)$  of  theorem \ref{exismes}}\label{infcons}

We first define a measure $\mu$ on ${\cal B}_{\infty}$. We consider
as before the set $B_{n}=B_{n}(v_{0})$, and pose $\mu_{n}=\mu^{B_{n}}$.

The set  ${\cal B}_{\infty}$ equipped with the product topology is
the projective limit of the sequence $({\cal W}(B_{n})$. We define $\mu$
as  the projective limit of the sequence $\{\mu_{n}\inn$. If $\mu^{3}$,
$\mu^{4}$ have full support on $O_{3}$ and $O_{4}$ respectively, then the measure  $\mu_{n}$ has full
support on  $\bar{\cal W}(B_{n})$. It follows that $\mu$ has full support on $\bar{\cal B}_{\infty}$.

The only non immediate  property of $\mu$ is the invariance under the ideal
triangle group $PSL(2,{\mathbb Z})$.

Notice first that if $g$ belongs to the stabiliser of the vertex $v_{0}$, we have that $g_{*}\mu_{n}=\mu_{n}$: this follows from the invariance of  $\mu^{3}$ under cyclic permutations, and from property $(iv)$ of proposition \ref{finitemes}.

Then, because  of the symmetries of $T$ and
the uniqueness of our construction, we have that if $g\in F$, $g_{*}\mu^{A,v}=\mu^{g(A),g(v)}$, and therefore
$g_{*}\mu^{A}=\mu^{g(A)}$, because of corollary \ref{fininv}.

It follows that $g_{*}\mu$ is the projective limit measure of the
projective limit of $\{{\cal W}(g(B_{n})\inn$, which is also ${\cal
B}_{\infty}$. 

To conclude, we just have to remark  that, thanks to $(v)$ of
proposition \ref{finitemes}, for whatever sequence of finite $v$-connected set
$\{D_{n}\inn$ in $B$, such that $D_{n+1}\subset D_{n}$ and $\cup_{n}D_{n}=B$,
the projective limit measure associated with the sequence of
$\{\mu^{D_{n}}\inn$ coincides with $\mu$.\bull

\subsection{Ergodicity}

We shall now prove property  $(vi)$ of theorem \ref{exismes}. Let's first
introduce some definitions. 
\subsubsection{Hyperbolic elements, pseudo-markov measure}
Let $F=PSL(2,{\mathbb Z})$ be the ideal triangle group, which we consider embedded in the isometry group
of the Poincar{\'e} disk. We  shall say $\gamma\in F$
is {\it hyperbolic}, if  $\g$ is an hyperbolic isometry. Notice that since
$F$ is Zariski dense, it contains many hyperbolic elements.

We also say a measure on ${\cal
B}^{0}_{\infty}/\Gamma$ is {\it pseudo-markov} if it satisfies the following property :
there exists an
integer $P$, such that for any $P$-disconnected and connected  subsets
$A$ and $C$ in $B$, if  $p_{A}$ and $p_{B}$ are the associated projections, then
$p_{*}^{A}\mu\otimes p_{*}^{B}\mu$  and  $(p^{A},p^{B})_{*}\mu$ are in the same measure class. By
proposition \ref{markov}, the measure we constructed in the last
section enjoys that property.

\subsubsection{Main result}
To conclude it suffices to prove:

\pro{Ergodicity}{Let $\mu$ be a $F$-invariant finite measure on 
 ${\cal
B}^{0}_{\infty}/\Gamma$, which is the pushforward of a pseudo-Markov
measure. Then $\mu$ is ergodic for the action of any hyperbolic
element of $F$, hence ergodic for $F$ itself.}
The proof is closely related to the proof of the ergodicity of subshifts of finite type, and is an avatar of Hopf's argument. We introduce {\it stable} and {\it unstable}  leaves in paragraph \ref{staunsta}, using {\it vanishing sequences} of sets defined in paragraph \ref{vani}. We finally conclude using Birkhoff Ergodic Theorem. 

\subsubsection{Hyperbolic elements of $F$}\label{vani}

Let $X$ be a topological space and $\gamma\in C^{0}(X,X)$, we shall say
for short that 
a sequence of non empty subsets $\{V_{n}\inn$  is a {\it vanishing
sequence} for $\gamma$ if:
\begitem
\iti  $V_{n+1}\subset V_{n}$;
\itii $\bigcap_{n\in {\mathbb N}} V_{n}=\emptyset$;
\itii for all compact subset $K$ of $X$, and $n\in {\mathbb N}$, there
exists $p\in {\mathbb N}$, such that $\gamma^{p} (K)\subset V_{n}$.
\enitem

\lem{hyperF}{Let  $\g$ be an  hyperbolic element in $F$. Then,  there exist two
families of connected subsets of 
$B$, $\{U^{+}_{n}\inn$ and $\{U^{-}_{n}\inn$,
, which are respectively vanishing sequences
for $\gamma$ and  for $\gamma^{-1}$, such that furthermore 
$U^{+}_{0}\cap U^{-}_{0}=\emptyset$.}

\preu  this is a consequence of elementary hyperbolic geometry. Indeed, if we
consider $F$ as a subgroup of the hyperbolic plane, the fixed
points of $\gamma$ on the boundary at infinity are not vertices of the
tiling by ideal triangles, and the lemma follows. \bull 

\subsubsection{Contractions}\label{staunsta}
Let now $\g$, $\{U^{\pm}_{n}\inn$ be as in  lemma \ref{hyperF}. We first introduce equivalence relations
amongst elements of $\bOi$. We say $f\sim^{+}_{n} g$, if $f\vert_{U_{n}^{+}}=g\vert_{U_{n}^{+}}$.
If $f\in {\cal B}^{0}_{\infty}$, let ${\cal F}^{+}_{n} (f)$ be the equivalence class of $f$.
Finally define $f\sim^{+} g$, if there exists $n$ such that  $f\sim^{+}_n
g$ and note ${\cal F}^{+}(f)$ the equivalence class of $f$. Observe that
$$
{\cal F}^{+}(f)=\bigcup_{n\in{\mathbb N}}{\cal F}^{+}_{n} (f),
$$
and let's define $\sim^-$, and ${\cal F}^{-}(f)$ in a symmetric way. These equivalence classes are going to play the role of the stable and unstable leaves of hyperbolic systems.

We shall prove:
\pro{contract}{There exist a $\G$-invariant metric on ${\cal
B}^{0}_{\infty}$ inducing the natural topology, such that for all
$f\in{\cal
B}^{0}_{\infty}$ and $g\in{\cal F}^{+}(f)$ then 
$$
\lim_{p\rightarrow+\infty} d(\g^{p}(f),\g^{p}(g))=0,
$$
and similarly if $g\in{\cal F}^{-}(f)$ then 
$$
\lim_{p\rightarrow+\infty} d(\g^{-p}(f),\g^{-p}(g))=0.
$$
}
\preu  
We first define a metric on ${\cal B}^{0}_{\infty}$ depending on the
choice of a vertex $v_{0}$ of the tree $T$. Let $B_{n}\subset B$ be
defined as in paragraph \ref{prooffinitemes}. Let ${\cal T}_{n}$ be the set
of tribones of $B_{n}$. Let $t$ be a tribone, then  
\begin{itemize}
\tit let ${\cal B}^{t}_{\infty}$ to be the set of
maps from $B$ to $W$, such that the image of $t$ lies in ${
O}^{3}$; 
\tit if $t\in {\cal T}_{n}$, let ${\cal B}^{t}_{n}$ be  the set of
maps from $B_{n}$ to $W$, such that the image of $t$ lies in ${
O}^{3}$; notice that $\Gamma$ acts properly on ${\cal B}^{t}_{n}$.
\end{itemize}
Next,
\begin{itemize}
\tit let $\delta_{n}^{t}$ be a $\G$-invariant distance of diameter less than 1 on  ${\cal
B}^{t}_{n}$ which induces the product topology;
\tit let $d_{n}^{t}$ the semi-distance on ${\cal B}^{t}_{\infty}$, induced from $\delta_{n}^{t}$ by the canonical projection; notice that  
the product topology of  ${\cal B}^{t}_{\infty}$ is induced by the
family of semi-distances $\{d_{n}^{t}\}_{n\in{\mathbb N},t\in{\cal T}_{n}}$.
\end{itemize}
By definition, recall that 
$$
{\cal B}^{0}_{\infty}=\bigcup_{\hbox{tribones}~t}{\cal B}^{t}_{\infty}.
$$
If $t\in{\cal T}_{n}$, we extend $d_{n}^{t}$ to  ${\cal B}^{0}$ in the following way: 
$$
\left\{
\begin{array}{cc}
d_{n}^{t}(x,y)&=0,~~\hbox{if}~x,y\in{\cal B}^{0}_{\infty}\setminus{\cal B}^{t}_{\infty},\\
d_{n}^{t}(x,y)&=1,~~\hbox{if}~y\notin{\cal B}^{t}_{\infty}, 
x\in{\cal B}^{t}_{\infty}.
\end{array}
\right.
$$
Ultimately, we define a $\G$-invariant metric $d$ on ${\cal B}^{0}_{\infty}$, by
the formula
$$
d(x,y)=\sum_{n\in{\mathbb n}}{\frac{1}{2^{n}}}\sum_{t\in{\cal
T}_{n}}{\frac{1}{\#{\cal T}_{n}}}d^{t}_{n}(x,y).
$$
By construction of this distance, it follows that is $f$ and $g$ coincides on $B_n$
then $d(f,g)\leq(\frac{1}{2})^{n-1}$. In particular,since 
$$
\forall q, \nin~, \exists p \in{\mathbb N}\hbox{~~such~that~~} \g_{p_n}(B_n)\subset U_q,
$$
it follows  that for every $n$, if $f\sim^{+}_q g$, then 
there exists $p\in{\mathbb N}$, such that 
$$
d(\g^{p}(f),\g^{p}(g))\leq (\frac{1}{2})^{n-1}.
$$ 
This 
ends the proof of the proposition .\bull 
\subsubsection{Preliminary steps for the proof of ergodicity}
Define for every
bounded function $\f$ on $\bOi$
$$
\f^{+}=\limsup_{n\rightarrow+\infty} (\f\circ\g^{n} ),
$$
and
$$
\f^{-}=\limsup_{n\rightarrow+\infty} (\f\circ\g^{-n} ).
$$
Let's first prove:
\pro{Anosov}{Let $\f$ be a continuous $\G$-invariant function on $\bOi$, such that the quotient function
on $\bOi /\G$ is compactly supported. Then $f\sim^+ g$ implies
$\f^+ (f)=\f^+ (g)$ and, $f\sim^- g$ implies $\f^- (f)=\f^- (g)$.}
\preu  Notice first that $\f$ is bounded and uniformly continuous. Hence, the proposition follows at once from
proposition \ref{contract}.\bull 

A second preliminary step is:
\pro{fubini}{Let $\m$ be a locally finite pseudo-Markov  measure of full support on $\bOi$. Let  $E$ be a set
of $\m$-full measure. Then for $\mu$-almost every $f$, there  exists a set $F_{f}$ of $\m$-full measure such that 
$$
\forall g \in F_{f},~~\exists h\in E,~~ s.t.~~f\sim^{+} h\sim^{-} g.
$$}
\preu  We should first notice that the set of equivalence classes of $\sim^{+}_n$ is precisely
${\cal W}(U^{+}_n )$, the space of maps from $U^{-}_n$ to $W$. A similar 
statement holds for  $\sim^{-}_n$. Fix some integer $n$, for which
$U^{+}_{n}$
and $U^{-}_{n}$ are $P$-disconnected. Let $p^+$ be
the natural continuous projection
$$
\bOi\mapsto {\cal W}(U^{+}_n ).
$$
Define $p^-$ a similar way. At last, let $p=(p^+ ,p^- )$.

If $E$ has full measure, then $p(E)$ has full measure for $p_{*}\m$. Hence by the pseudo-Markov property, it has full measure for $p^{+}_{*}\m\otimes
p^{-}_{*}\m$.

From Fubini Theorem, we deduce there is a set of full measure $A$ in
${\cal W}(U^{+}_n )$, such that for every $a\in A$, the set
$$
V_{a}=\{c\in {\cal W}(U^{-}_n ),~(a,c)\in p(E)\}
$$
has full measure for $p^{-}_{*}\m$.

In particular, for every $f\in (p^{+})^{-1}(A)$, the set
$F_{f}=(p^{-})^{-1}(V_{p^{+}(f)})$ has full measure.

Now, by construction if $f\in (p^{+})^{-1}(A) $ and $g\in F_{f}$, then $p^{-}(g)=p^{-}(h)$, where $h\in
E$ and $p^{+}(h)=p^{+}(f)$. This is exactly what we wanted to
prove.\bull 

\subsubsection{End of the proof of ergodicity}
In this paragraph, we will prove proposition \ref{Ergodicity}.
Let $\g$ be some hyperbolic element in $F$. Let $\m$ be the $F$-invariant measure on  $\bOi /\G$, constructed in \ref{exismes}. We first  use the
ergodic decomposition theorem and write
$$
\m=\int_{Z}\n_{z}d\l (z),
$$
where for $\l$-almost every $z$ in $Z$, $\n_{z}$ is  an ergodic
measure for $\g$.

To conclude, it suffices to show that for any  continuous and compactly
supported function $\psi$ on $\bOi/\G$, and for every $z$ and $u$ in
$Z$, we have $\int \psi dv_{z}=\int \psi dv_{u}$.

Let now $\f=\psi\circ \pi$, where $\pi$ is the natural projection
from $\bOi$ to $\bOi/\G$. Let's define as for proposition
\ref{Anosov}, the measurable functions $\f^{+}$ and $\f^{-}$.
From Birkhoff ergodic theorem, we deduce that for $\n_{z}$-almost
every $x$,  if  $\pi(y)=x$, we have 
$$
\f^{+}(y)=\f^{-}(y)=\int_{\bOi/\G}\psi d\n_{z}.\leqno{(*)}
$$
In particular, there exists a set of $\m$-full measure $E$ on which  $\f^{+}=\f^{-}$. 

Now, we apply proposition \ref{fubini}, and deduce that for $\m$-almost
every $x$, there exists a set of full measure $F_{x}$ with the following
property: if $b\in
F_{x}$ then  there exists $a\in E$ such that $x\sim^{+} a\sim^{-} b$.

From proposition \ref{Anosov}, we deduce that $\f^{+}(a)=\f^{+}(x)$
and $\f^{-}(b)=\f^{-}(a)$. From the definition of $E$, we get that
$\f^{-}$ is constant and equal to $\f^{-}$ on $F_{x}$, hence $\m$-almost everywhere.

Using $(*)$, we ultimately get that for almost every $z,~u\in Z$
$$
\int_{\bOi/\G}\psi d\n_{z}=\int_{\bOi/\G}\psi d\n_{u},
$$
which is what we wanted to prove.\bull

\section{Configuration data and the boundary at infinity of a hyperbolic 3-manifold}
\label{crossratio}
We describe here our main, and actually unique useful example : the Markov configuration data associated 
to a hyperbolic 3-manifold.

Let in general $\di$ be the boundary at infinity of a negatively curved
3-manifold $M$.  Let $\Gamma$ be  a discrete, torsion free and cocompact group of isometries of $M$.

Unless otherwise specified, we shall assume $M$ is the hyperbolic 3-space $\yp$. Then,  $\di=\dih$ is canonically identified with $\cpi$. In this identification, the action of the group of isometries of $M$ on $\di$ coincides  with the action of $PSL(2,{\mathbb C})$ on $\cpi$.

As we explained in \ref{mainexample}, the (3,4)-configuration data we shall study is the following :
\begin{itemize}
\tit $W=\di=\cpi$,
\tit $O_{3}$ is the subset of
$\di^{3}$ consisting of triples of different points :
$$
O_{3}=\{(x,y,z)\in \di~ /~x\not=z\not=y\not=x\}.
$$
\tit  $O_{4}$  is  the set of points whose crossratios have a non zero
imaginary part;
\end{itemize}

We have
\pro{datahyp}{The quadruple
$(\cpi,\G, O_{3},O_{4})$ is a Markov (3,4)-con\-fi\-gu\-ra\-tion data.} 

It is obvious. The only point that requires a check is hypothesis (e). In the last paragraph \ref{crossratio2}, we will devise a fancy (and far too long) proof of this fact. Of course, a straightforward check would give that this configuration data satisfies (e) for $N=10$, instead of $N=1000$, provided by our proof. However, I hope the scheme of this proof might be useful in more general situations.

In the next paragraph we explain how to turn this example in a regular measured configuration data in many ways, using {\it equivariant family of measures} (cf \ref{equifamili}).

\subsection{Measured configuration data}
\label{mesconfdata}
In view of \ref{datacheck}$(ii)$ we need  to produce $\mb^{3}$ in the
Lebesgue class of $m\otimes m\otimes m$ for some measure class $m$ of
full support, such  that the pushforward of $\mb^{3}$ on
$O_{3}/\G$ is finite. Then we have to build a $\G$-equivariant map $\nb$: 
\begin{displaymath}
\left\{ \begin{array}{cc}
O_{3} & \rightarrow {\cal P}_{m}(W) \\
(a,b,c) & \mapsto\nb_{(a,b,c)} 
\end{array}
\right.
\end{displaymath}
where ${\cal P}_{m}(W)$ is the set of finite radon measures on
$W$ in the measure class of $m$.

We shall do that using the notion of equivariant family of measures
described  by F. Ledrappier in \cite{led}, and which is
a generalisation of that of conformal densities due to 
D. Sullivan in \cite{Su}.

\subsubsection{Equivariant family of measures}\label{equifamili}
An {\it equivariant family of measures on the boundary} a map $\mu$ which
associates to every $x\in M$ a finite measure $\mu_{x}$ on $\di$
such that: 
\begitem
\iti for all $\g$ in $\G$, $\m_{\g x}=\g_{*} \m_{x}$,
\itii For all $x,y\in M$, $\m_{x}$ and $\m_{y}$ are in the same Lebesgue
class.
\enitem
In particular we can write
$d\m_{x}(a)=e^{-\g_{a}(x,y)}d\m_{y}(a )$.
Actually, the original definition requires some regularity of the
function $(a , x,y)\mapsto \g_{a}(x,y)$, which we shall not need in
the sequel.

A typical example arises when one associates to a point $x$ the pushforward by
the exponential map of the Liouville measure on the unit sphere at $x$. 

When $c_{\eta}(x,y)=\delta B_{\eta }(x,y)$, where $B_{\eta}(x,y)$ is the
the {\it Busemann function} defined by
$$
B_{\eta} (x,y)=\lim_{z\rightarrow \eta}(d(x,z)-d(y,z))
$$
the corresponding equivariant family of measures is called a {\it conformal
density of ratio $\delta$}. Among these  is the  Paterson-Sullivan
measure.

In \cite{led} which also contains many references to related results,
F. Ledrappier discusses various ways of building
equivariant families of measures, and in particular relates them to
other notions
like {\it crossratios, Gibbs currents, transverse invariant measures to the
horospherical foliations  etc}. As a conclusion, there exist
numerous
examples of  equivariant families
of measures, all mutually singular.
\subsubsection{End of the construction}
Let's go back to our construction now. Let first $\beta$
$$
(a,b,c )\mapsto \beta_{a,b,c}
$$
be a
$\G$-equivariant map from $O_{3}$ to $M$. For instance, we can take the
barycentre of the sum of the three Dirac measures concentrated at $a$,
$b$, and $c$. Now define, if $x\in M$

$$
d\mb^{3}(a,b,c)=e^{c_{a}(x,\beta_{a,b,c})+c_{b}(x,\beta_{a,b,c})+c_{c}(x,\beta_{a,b,c})}d\mu_{x}\otimes d\mu_{x}\otimes d\mu_{x} (a,b,c).
$$

It follows from the definition of equivariant families of measure  that this
definition is independent on $x$. If $\Gamma$ is a group of isometries
then $\mb^{3}$ is $\Gamma$-invariant. Furthermore, if $\Gamma$ is cocompact then  the
corresponding measure is finite on $O_{3}/\Gamma$.

For $\nb$, we can now just take the map $(a,b,c )\mapsto \mu_{\beta_{a,b,c}}$.

\subsubsection{Negatively curved  3-manifolds}
We do not have used previously the hyperbolic structure. Let's take for a general negatively curved $M$ and cocompact group of isometries $\Gamma$
\begin{itemize}
\tit $W=\di$,
\tit $O_{3}=\di^{3}\setminus\Delta_{3}$,
\tit $O_{4}$ any $\lambda_{4}$-invariant subset such that $p(O_{4})=O_{3}$.
\end{itemize}
Then the  previous construction works provided that $O_{4}$ is invariant under  all $\s_{4}^{+}$.  For instance, we could take ${O}_{4}=U_{4}$, but the corresponding construction seems to be of no use for our problem. 
For the moment, I have not been able to construct a configuration data  adapted to the problem, for general negatively curved 3-manifolds.

\subsection{Complex crossratio}
\label{crossratio2}
Let  $[a;b;c;d]$ be the
complex crossratio of four points of ${\mathbb C\mathbb P}^1$, such that $[0;1;\infty;z]=z$.  Let $\Im(\a)$ be the imaginary part of the complex number $\a$.

We will single out the geometric properties of the crossratio which are actually used in proposition \ref{datahyp}.
 
\subsubsection{disks}
We  associate to every triple $(a,b,c)$, a {\it disk}  $D(a,b,c)$, defined  by
$$
D(a,b,c)=\{z \in {\mathbb C\mathbb P}^1 ~~/~~ {\Im}( [a;b;c;z] ) < 0\}.
$$ 

If we consider $D$  as  a map  from $O_{3}$ to the set of  subsets of $\di$,
it enjoys the following properties:
\smallskip
\begitem\label{virdef}
\iti the map $D$ is  $\G$-equivariant;
\itii $D(a,b,c)\cup D(a,c,b)$ is dense;
\itiii for every $(a,b,c)$ in $O_{3}$, $a$ belongs to the closure of $D(a,b,c)$;
\itiv  Let $O_{4}=\{a,b,c,d~ /~(a,b,c)\in O_{3}~,~~d\in
D(a,b,c)\}$, then $O_{4}$ is an open set invariant under 
the oriented permutations $(a,b,c,d)\mapsto (d,c,b,d)$.
\enitem
In our particular case, all these properties are easily checked from the invariance of the
crossratio and it's behaviour under permutations.

We will prove:

\pro{virdata}{Let  $D$  be satisfying properties $(i)$ to $(iv)$   of \ref{virdef}. Then the quadruple
$(\di,\G, O_{3},O_{4})$ is a Markov (3,4)-configuration data.} 

In our very precise situation, we could devise a quick proof of that fact. However, we will give a somewhat longer proof: the idea is to stress the importance of properties  $(i)$ to $(v)$   of \ref{virdef}, and forget a while the complex structure on $\di$.

{\it Proof of the proposition:}   it only remains to prove (e) of
definition \ref{data} which characterises Markov configuration data. Let's recall it:
\begitem
\ite  There exist some constant $p\in {\mathbb N}$, such that if $(a,b,c)$ and
$(d,e,f)$ both belong to $O_{3}$, then there exists a sequence
 $(q_{1},\ldots,q_{j})$ of elements  of $O_{4}$, where $j\leq p$, such
 that if
 $q_{n}=(q^{1}_{n},q^{2}_{n},q^{3}_{n},q^{4}_{n})$ 
 then $(q_{1}^{1},q_{1}^{2},q_{1}^{3})=(a,b,c)$,
 $(q_{j}^{3},q_{j}^{2},q_{j}^{4})=(d,e,f)$ and at last
 $(q^{2}_{n},q^{4}_{n},q^{3}_{n})=(q^{1}_{n+1},q^{2}_{n+1},q^{3}_{n+1})$
 or $(q^{3}_{n},q^{2}_{n},q^{4}_{n})=(q^{1}_{n+1},q^{2}_{n+1},q^{3}_{n+1})$ ,
\enitem
In order to proceed towards a proof we shall write that $(a,b,c)
\stackrel{p}{\leadsto}(d,e,f) $ if $(a,b,c)$ and
$(d,e,f)$ satisfies this condition (e).
With these  notations at hands, one immediately checks:
\begitem
\item[-]
$(a,b,c)\stackrel{p}{\leadsto}(u,v,w)$ implies $(w,u,v)\stackrel{p}{\leadsto}(b,c,a)$;
\item[-] {\it composition rule}:
if
$t_{1}\stackrel{p}{\leadsto} (a,b,c)$, and $(a,b,c)\stackrel{q}{\leadsto} t_{3}$ or $(b,c,a)\stackrel{q}{\leadsto} t_{3}$ then
$t_{1}\stackrel{p+q}{\leadsto} t_{3}$;
\item[-] 
$(a,b,c,d)\in O_{4}$ exactly means that $(a,b,c)\stackrel{1}{\leadsto}(c,b,d)$.
\enitem
We are going to proceed through various steps. 

\medskip
\noindent
{\bf Step 1:} {\em for any $(a,b,c)$ there exists $(a_{1},b_{1},c_{1})$
arbitrarily close to $(a,b,c)$ such that 
$(a,b,c)\stackrel{3}{\leadsto}(a_{1},c_{1},b_{1})$.}

\medskip
We shall prove that using property $(iii)$ of the definition
\ref{virdef}. First, using $(iii)$ of \ref{virdef}, we choose $b_{1}$ arbitrarily close to $b$ such
that $(b,c,a,b_{1})\in O_{4}$. 

Next, using $(iii)$ again, we choose $c_{1}$ arbitrarily close
to $c$ such that $(c,b,a,c_{1})\in O_{4}$, and still, because $O_{4}$ is
open,  $(b,c_{1},a,b_{1})\in O_{4}$. 

At last, using $(iii)$ again, we
choose $a_{1}$ arbitrarily close to $a$ such that $(a,b,c,a_{1})\in
O_{4}$ and still $(b,c_{1},a_{1},b_{1})$ and $(c,b,a_{1},c_{1})$ in
$O_{4}$. 

It follows that we have 
\begin{eqnarray*}
(a,b,c)&\leadsto&(c,b,a_{1})\\
(c,b,a_{1})&\leadsto&(a_{1},b,c_{1})\\
(b,c_{1},a_{1})&\leadsto&(a_{1},c_{1},b_{1}).
\end{eqnarray*}
The composition rule implies now that $(a,b,c)\stackrel{3}{\leadsto}(a_{1},c_{1},b_{1})$.

\medskip
\noindent
{\bf Step 2:} {\em for any $(a,b,c)$ there exists
$(a_{1},b_{1},c_{1})$ arbitrarily close to $(a,b,c)$ such that
$(a,b,c)\stackrel{3}{\leadsto}(b_{1},a_{1},c_{1})$.}

\medskip
The proof is symmetric: first we notice,  using $(iii)$, that we can find $c_{1}$
arbitrarily close to $c$ such
that 
\begin{eqnarray*}
(c,a,b)&{\leadsto}&(b,a,c_{1}).
\end{eqnarray*}

We choose $b_{1}$, arbitrarily close to $b$, such that
\begin{eqnarray*}
(b,a,c)&{\leadsto}&(c,a,b_{1})\\
(c,a,b_{1})&{\leadsto}&(b_{1},a,c_{1}).
\end{eqnarray*}
Lastly, we choose $a_{1}$, arbitrarily close to $a$, such that
\begin{eqnarray*}
(a,b,c)&{\leadsto}&(c,b,a_{1})\\
(b,a_{1},c)&{\leadsto}&(c,a_{1},b_{1})\\
(c,a_{1},b_{1})&{\leadsto}&(b_{1},a_{1},c_{1}).
\end{eqnarray*}
The composition rule implies the desired statement.

\medskip
\noindent
{\bf Step 3:} {\em for any $(a_{1},b_{1},c_{1})$ close enough to
$(a,b,c)$ then $(a,b,c)\stackrel{36}{\leadsto}(a_{1},b_{1},c_{1})$.}

\medskip
It suffices to prove that $(a,b,c)\stackrel{36}{\leadsto}(a,b,c)$.
From the first step, we have that given $(a,b,c)$  there exists $(a_{4},b_{4},c_{4})$
arbitrarily close to $(a,b,c)$ such that 
$$
(a,b,c)\stackrel{3}{\leadsto}(a_{4},c_{4},b_{4}).
$$
Next applying the first step one more time, we can choose $(a_{3},b_{3},c_{3})$ arbitrarily close to $(a_{4},b_{4},c_{4})$
such that
$$
(a_{4},c_{4},b_{4})\stackrel{3}{\leadsto}(a_{3},b_{3},c_{3}),
$$ 
and this
leads to 
$$
(a,b,c)\stackrel{6}{\leadsto}(a_{3},b_{3},c_{3}).
$$
Actually
we can choose  $(a_{3},b_{3},c_{3})$ close enough  to
$(a_{4},b_{4},c_{4})$ such that still
$$
(a,b,c)\stackrel{3}{\leadsto}(a_{3},c_{3},b_{3}).
$$
At last, we choose thanks to step 2, $(a_{2},b_{2},c_{2})$ arbitrarily
close 
to $(a_{3},b_{3},c_{3})$
such that
$$
(a_{3},c_{3},b_{3})\stackrel{3}{\leadsto}(c_{2},a_{2},b_{2}),
$$ 
and this
implies
$$
(a,b,c)\stackrel{6}{\leadsto}(c_{2},a_{2},b_{2}),
$$
As before,  we can choose  $(a_{2},b_{2},c_{2})$ close enough so that we still
have
$$
(a,b,c)\stackrel{6}{\leadsto}(a_{2},b_{2},c_{2}).
$$
From this last  relation, we obtain
$$
(c_{2},a_{2},b_{2})\stackrel{6}{\leadsto}(b,c,a).
$$
It follows
$$
(a,b,c)\stackrel{12}{\leadsto}(b,c,a),
$$
Hence 
$$
(a,b,c)\stackrel{36}{\leadsto}(a,b,c).
$$

\medskip 
\noindent
{\bf Step 4:} {\em for any $a,b,c$ and any permutation $\sigma$ we have
$$
(a,b,c)\stackrel{100}{\leadsto}(\sigma(a),\sigma(b),\sigma(c)).$$}

\medskip 
This follows easily from the previous steps.

\medskip 
\noindent
{\bf Step 5:} {\em for any $a,b,c$ there exist an open dense set of
$d$ such that 
$$
(a,b,c)\stackrel{300}{\leadsto}(b,c,d).
$$}

\medskip 
Indeed, from hypothesis $(ii)$ of \ref{virdef}, there exist an open dense
set of $d$ such that either $(a,b,c)\stackrel{1}{\leadsto}(c,b,d)$ - in which case
we are done - or
$(a,c,b)\stackrel{1}{\leadsto}(b,c,d)$ and we obtain our assertion using step 4
twice.

\medskip 
\noindent
{\bf Final step:} {\em for any $(a,b,c,d,e,f)$, we have $(a,b,c)\stackrel{1000}{\leadsto}(d,e,f)$.} 
\medskip
Using step 5 three times,  we obtain there exists an open dense set of
$(u,v,w)$ such that $(a,b,c)\stackrel{900}{\leadsto}(u,v,w)$, hence our
conclusion thanks to step 3.

\bigskip
The proof is complete although, obviously, 1000 is not the optimal
constant. Also this proof is far too complicated in our case, but
one of my  hope is to build a map $D$ satisfying $(i)$, $(ii)$,  $(iii)$
$(iv)$ and $(v)$ of \ref{virdef} for a general negatively curved 3-manifold.\bull

\sk{Convex surfaces and configuration data}\label{coding}

Let's again $N=M/\G$ be a compact hyperbolic 3-manifold. In the last
section we have built a configuration data associated to that situation,
and we can extend that to a measured configuration data in many ways
({\it cf} \ref{mainexample})

We consider now the restricted configuration space $\bbi$, associated
to  that situation with it's measure $\mu$, invariant under $\G$ and ergodic
under the action of the ideal triangle group $F$ as constructed in
\ref{exismes}.

We turn $\bbi$ into an ergodic  
Riemannian lamination by a suspension procedure,
namely we consider
$$
{\cal F}=(\bbi\times \hyp )/F,
$$
where $F$ acts as an isometry group on $\hyp$ and diagonally on
$\bbi\times\hyp$.

The ergodic and $\G\times F$-invariant measure $\mu$ gives rise to a
transversal  $\G$-invariant and ergodic measure on ${\cal F}$ that we
shall also call $\mu$.  

Our aim is now to prove:

\pro{exisphi}{There exist a continuous leaf preserving map 
$\Phi$ with dense image from $\cal F/\G$ to ${\cal N}$, the space of
$k$-surfaces in $N$.}

This proposition, loosely speaking, explains that our combinatorial
construction codes for convex surfaces. As a corollary, we obtain our
main theorem

\theo{maintheo}{Let $N=M/\G$ be a compact negatively curved 3-manifold whose metric can be deformed through negatively curved metrics to a hyperbolic one. Then  there
there exist infinitely many  mutually singular ergodic transversal measures of full
support on ${\cal N}$, the space of $k$-surfaces of $N$.}

\subsection{Bent and pleated surfaces}\label{bout}

Recall that a {\it $\cpi$-surface} is a surface locally modelled on $\cpi$.

We  shall have to recall facts about (locally convex) pleated surfaces, mostly without demonstrations,  especially when dealing with  the relation between measured geodesic laminations and $\cpi$-structures which has been described by W. Thur\-ston. A useful reference is \cite{tan}, where H.~Tanigawa  gives a description and some results about this relation.

The main fact about this construction is the following:  to every hyperbolic surface $S$ (maybe non complete) and every measured geodesic lamination $\mu$ we can associate 
\begin{itemize}
\item a pleated locally convex surface in the hyperbolic space, 
\item a 3-manifold $B(S,\mu)$, the {\it end} of $(S,\mu)$, 
\item a $\cpi$-surface $\Sigma$ which is going to be the boundary at infinity of the end. 
\end{itemize}
The map which associates to $(S,\mu)$ the $\cpi$-surface $\Sigma$ is called the {\it Thurston map}, we shall denote it by $\Theta$. Notice that since $S$ is not assumed to be complete, this map has no reason to be injective.

\sssk{An example} For the sake of completeness, we briefly recall the Thurston's construction in a special case, which will be the one we shall actually need. 

Let $S$ be a open subset of $\hyp$  (maybe non complete) which is the union of totally geodesic ideal polygons. To every edge $e$ of this tiling, we associate a positive number $\theta_{e}$ less than $\pi$. The data $\mu$ consisting of the edges of the tiling and of the assigned positive numbers is  a specific example of a geodesic lamination.

We may think of every polygon $T$ as totally geodesically embedded in $\yp$. Let $n_{T}$ be the exterior normal field along $T$. Let $p_{T}$ be the map from $T\times]0,\infty[$ defined by 
$$
p_{T}~:~(x,s)\mapsto \exp(sn_{T}(x)).
$$
Let $P_{T}$ be the {\it prism} over $T$, {\it i.e.} the image of $p_{T}$. 

Let $e$ be an edge of the tiling of $S$, intersection of two polygons $T^{e}_{0}$ and $T^{e}_{1}$ and considered as a geodesic in $\yp$. The $\theta_{e}$-{\it wedge} over $e$ is the closed set delimited by the two half-planes whose boundary is $e$ and forming an angle $\theta_{e}$.

Finally, the {\it end} $B(S,\mu)$ of $({S},\mu)$ is the reunion of all prisms and edges. Notice that there is a canonical isometric local homeomorphism from $M_{\Sigma}$ to $\yp$. 

In this special case, $S$ is isometrically immersed in $\yp$ as a pleated surface.

\sssk{Facts}

The following propositions, whose proof follows from  results explained in 
\cite{tan}, summarises the property of Thurston's construction we shall need in the sequel. All these properties rely on the following observation:

\obser{}{Let $S$ be a hyperbolic surface (maybe non complete) and $\mu$ a geodesic lamination, let $D$ be an embedded $\cpi$-disk in $\Theta(S,\mu)$, then there is an embedding of the (hyperbolic) half space $P$  in $B(S,\mu)$ such that, the boundary at infinity of this embedded $P$ is precisely $D$.}

From this observation, we deduce easily the following results.

\pro{thur1}{Let $S$ be a (maybe non complete) hyperbolic surface and $\mu_{1}$ a measured geodesic lamination on $S$. Let $\mu_{2}$ be a measured geodesic lamination on $\hyp$ supported on finitely many geodesics. Assume  $\Theta(\hyp,\mu_{2})$ injects by $f$ (as a $\cpi$-surface) in $\Theta(S,\mu_{1})$. Then $B(\hyp,\mu_{2})$ injects in $B(S,\mu_{1})$ in such a way the associated injection between the boundary at infinity is $f$.}

\pro{thur2}{Let $M$ be a $\cpi$-surface. Then there exists an exhaustion of $M$ by relatively compact $\cpi$-surfaces $M_{i}$ such that $M_{i}=\Theta(\hyp,\mu_{i})$ where $\mu_{i}$ is supported on finitely many geodesics.}

Let $S$ be a locally convex immersed surface in $\yp$. Let $n$ be its exterior normal field. We define the {\it end}, to be the 3-manifold  $B_{S}$ diffeomorphic to $S\times]0,\infty[$ equipped with the hyperbolic metric induced by  the immersion  $(s,t)\mapsto exp(tn(s))$. 
In particular, the boundary at infinity $S_{\infty}$ of $B_{S}$ is a 
$\cpi$-surface.

\pro{thur3}{Let $S$ be a locally convex immersed surface in $\yp$ such that $B_{\infty}$ injects as a $\cpi$-surface in $\Theta(\hyp,\mu)$. Then $B_{S}$ injects in $B(\hyp,\mu)$.}

\subsection{Tilings and related definitions}

We shall denote by $T(a,b,c)$ the ideal triangle in $\hyp$ whose vertices are 
$a$, $b$ and $c$ in $\partial_{\infty}\hyp=\mathbb{RP}^1$.

Recall that we consider $\hyp$ periodically tiled by ideal triangles. Let's denote ${\cal T}^{0}$ the collection of ideal triangles of this triangulation. 
The set $B=\mathbb{QP}^{1}$ is the set of vertices at infinity of this triangulation ({\it cf} \ref{3tree}). 

Notice now that every monotone map $g$ from $B=\mathbb{QP}^{1}$ to $\partial_{\infty}\hyp=\mathbb{RP}^1$ defines a tiling by ideal triangles of an open set $U_{g}$ of $\hyp$. This triangulation is given by the collection ${\cal T}^{g}$ of triangles defined by
$$
{\cal T}^{g}=\{T(g(a),g(b),g(c))/T(a,b,c)\in{\cal T}^{0}\}.
$$

With these notations we have: 
$$
U_{g}=\bigcup_{T\in{\cal T}^{g}}T.
$$

\sssk{Pleated surfaces and tilings}

We will prove the following immediate  proposition: 

\pro{pleated}{For every $f\in \bbi$, there exists a unique monotone map $g(f)$ from $B$ to  $\partial_{\infty}\hyp$, a unique map $\psi_{f}$ from
 $U_{g(f)}$ to $\yp$, such that its restriction to  every tile $T(a,b,c)$ is totally geodesic, and the ideal triangle $\psi_{f}(T(a,b,c))$ has $f(a)$, $f(b)$ and $f(c)$ as vertices at infinity. Furthermore, $S^{0}_{f}=\psi_{f}$ is locally convex and $\psi_{f}$ depends continuously on $f$.}

Using this proposition, we introduce the following notations: 
we shall denote by $\mu_{f}$ the geodesic lamination on $U_{g(f)}$ whose support is the set of edges of  ${\cal T}^{g}$, each edge being labelled with the angle of the two corresponding triangles in $\yp$. Then we shall write $B_{f}$ for $B(U_{g(f)},\mu_{f})$ and $S_{f}^\infty$ for $\Theta(U_{g(f)},\mu_{f})$.

\preu The construction of $\psi_{f}$ is described in the statement of the proposition. The only point to check is the local convexity of $S^{0}_{f}$. This follows at once from the following observation: three points $(a,b,c)$ at infinity in $\yp$ determine an oriented totally geodesic plane $P$ in $\yp$, and the points in $\cpi$ ``below'' $P$ are precisely those points $d$ such that $\Im(a,b,c,d)<0$. \bull


\sssk{Tiling map}

Later on, we shall need  a technical device, called a {\it tiling map}, associated to every element of $\bbi$.

Let $C^{0}(\hyp,\cpi)$ be the space of continuous maps from $\hyp$ to $\cpi$ 
 with the topology of uniform convergence on every compact set. We use the notations of the previous section \ref{bout}.

The following proposition is obvious.

\pro{tiling_map}{There exists a continuous map $\xi$ 
\begin{eqnarray*}
\left\{
\begin{array}{rcl}
\bbi&\rightarrow&C^{0}(\hyp,\cpi)\\
f&\mapsto&\xi_{f}.
\end{array}
\right.
\end{eqnarray*}
which satisfies the following properties 

\begitem
\iti there exists an homeomorphism $h_{f}$ from $\hyp$ to $S^{\infty}_{f}$, such that  $\xi_{f}=i_{f}\circ h_{f}$,
\itii for every $T(a_{1},a_{2},a_{3})$ in  ${\cal T}^{0}$, the map  $\xi_{f}\vert_{T(a_{1},a_{2},a_{3})}$
extends continuously to $\{a_{1},a_{2},a_{3}\}$ in such a way that $\xi_{f}(a_{i})=f(a_{i})$, 
\itiii for every element $\gamma$ in $F$, $\xi_{f\circ\gamma}=\xi_{f}\circ\gamma$.
\enitem
}

By definition, $\xi_{f}$ is a {\it tiling map} associated  to $f$.

\ssk{$k$-surfaces and asymptotic Plateau problems}\label{asym}

We recall definitions and results from \cite{la} that we specialise in the case of $\yp$. 

Let $S$ be a locally convex surface immersed in $\yp$. Let $\nu_{S}$ be the exterior normal vector field to $S$. The {\it Gauss-Minkowski} (Figure \ref{minko}) map from  $S$ to 
$\partial_{\infty}\yp$ is the local homeomorphism $n_{S}$:
\begin{eqnarray*}
\left\{
\begin{array}{rcl}
S&\rightarrow&\partial_{\infty}\yp\\
x&\mapsto&n_{S}(x)=\exp(\infty \nu_{S}(x)).
\end{array}
\right.
\end{eqnarray*}   
\begin{figure}[h!]  
  \centerline{\psfig{figure=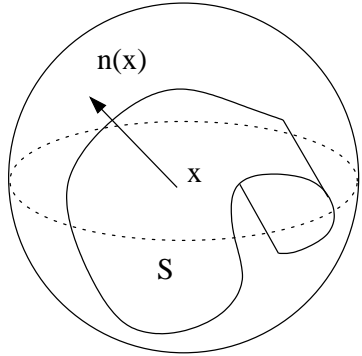}}
  \caption{\it Gauss-Minkowski map}\label{minko}
\end{figure}

An {\it  asymptotic Plateau problem} is a pair $(i,U)$ where $U$ is a surface, and $i$ is  a local homeomorphism from $U$ to $\partial_{\infty}\yp$.  A {$k$-solution} to an asymptotic Plateau problem $(i,U)$, is a $k$-surface $S$ immersed in $\yp$, such that there exists a homeomorphism $g$ from $U$ to $S$ such that $i=n_{s}\circ g$.

We proved (theorem A of \cite{la}) that there exists at most one solution of a given asymptotic Plateau problem. We also proved (theorem E of \cite{la}) that if $(i,U)$ is an asymptotic Plateau problem, and if $O$ is a relatively compact open set of $U$, then $(i,O)$ admits a solution.

We need the following proposition which uses the notations of section \ref{bout}

\pro{exisplat}{Let $f$ be an element of $\bbi$, then the asymptotic Plateau problem $(i_{f},S^{\infty}_{f})$ admits a $k$-solution.}

\preu Using proposition \ref{thur2}, we have an exhaustion of $S_{f}^{\infty}$, by $\cpi$-surfaces $M_{i}$, such that $M_{i}=\Theta(\hyp,\mu_{i})$ where $\mu_{i}$ is supported on finitely many geodesics. Let $B_{i}=B(\hyp,\mu_{i})$. According to \ref{thur1}, we have 
$$
B_{i}\subset B_{i+1}\subset B_{f}.
$$

Since $M_{i}$ is relatively compact in $S_{f}^{\infty}$, there exists, according to theorem E of \cite{la}, a $k$-solution $\Sigma_{i}$ to the asymptotic Plateau problem defined by $M_{i}$. According to proposition  \ref{thur3}, 
$$
B_{\Sigma_{i}}\subset B_{i}\subset B_{f}.
$$
Let now
$$
W=\bigcup_{i\in\mathbb N} B_{\Sigma_{i}}.
$$
We wish now to prove that $\partial W$  the boundary of $W$ is a $k$-surface solution of the asymptotic Plateau problem defined by $S_{f}^{\infty}$.

First we should notice that since $B_{\Sigma_{i}}\subset B_{f}$, there exists 
a constant $A$ just depending on $f$, such that every ball $\Sigma_{i}(x,A)$ of centre $x$ and radius $A$ in $\Sigma_{i}$, when considered immersed in $\yp$, is a subset of the boundary of a convex set. 
Hence, according to lemma 5.4.($iii$) of \cite{la2}, there exists a constant $C$ such that  we have, if $H$ is the mean curvature of $\Sigma_{i}$ and $d\sigma$ the area element: 
\begin{displaymath}
\int_{\Sigma_{i}(x,A)}H d\sigma\leq C.
\end{displaymath}

Let now $x_{i}$ be  a point in $\Sigma_{i}$. Assume the sequence  $\{x_{i}\}_{i\in \mathbb N}$ converges in the metric completion of $B_{f}$ to a point $x_{0}$. We conclude from theorem D of \cite{la2} that $\{(\Sigma_{i},x_{i})\}_{i\in\mathbb N}$ converges smoothly to a pointed $k$-surface $(\Sigma_{\infty},x_{0})$. 
It follows that $x_{0}$ is in the interior of $B_{f}$. Indeed, let $\partial B_{f}$ be the boundary of the metric completion of $B_{f}$. Notice that every point of $\partial B_{f}$ is included in an open geodesic segment drawn on $\partial B_{f}$. If $x_{0}$ belongs to  $\partial B_{f}$, it would  follow that the corresponding open segment is actually drawn on $\Sigma_{\infty}$ and this is impossible.

This argument finally shows that $\partial W$ is a $k$-surface, and that every half infinite geodesic joining a point of $\partial B_{f}$ to a point of $S_{f}^{\infty}$ intersects $\partial W$.
The conclusion follows. \bull

\ssk{Construction of the map $\Phi$}

In this section, we summarise the previous sections and  build a continuous map $\Phi$ from ${\cal F}$ to ${\cal N}$.

Let $f\in\bbi$. Let $\xi_{f}$ the tiling map of $f$ ( cf \ref{tiling_map} ).
Let $\Sigma_{f}$ be the $k$-solution of the asymptotic Plateau problem $(i_{f},S^{\infty}_{f})$ ( {\it cf} \ref{exisplat} ). Let $n_{f}$ be the Gauss-Minkowski map of $\Sigma_{f}$. We define $\Phi$ by
$$
\Phi([f,x])=(\Sigma_{f},n_{f}^{-1}(\xi_{f}(x))).
$$
Continuity follows from the uniqueness of the solution of an asymptotic Plateau problem.

\ssk{Density of the image of $\Phi$} 

The only point left to be proved in proposition \ref{exisphi} is the density of the image of $\Phi$. 

We start with an observation.
Let $S$ be a compact surface, ${\tilde S}$ its universal cover. Let 
 $\mu_{1}$ be a measured lamination on $S$ supported on finitely many geodesics. Assume the weight of every geodesic is strictly less than $\pi$.  Then, from the construction explained above we deduce that  $\Theta({\tilde S},\mu_{1})$ lies in the image of $\Phi$.

According to \ref{ksur}, the union of compact leaves of ${\cal N}$ is dense. It therefore suffices to prove that every compact $k$-surface belongs to the closure of the image of $\Phi$.

Let $S$ be such a compact $k$-surface in $N$. The underlying surface admits a $\cpi$-structure induced by the Minkowski-Gauss map. According to Thurston's parametrisation theorem \cite{tan}, such a surface is of the form $\Theta(S,\mu_{0})$ for a certain measured lamination $\mu_{0}$. We proved, using different words  (Corollary 1 of \cite{la3}) that the map which associates to every $\cpi$-structure on a compact surface, the $k$-surface solution of the corresponding asymptotic Plateau problem, is continuous. To complete our proof,  we just have to remark that the set of measured geodesic laminations with finite support  and such that the weight of every geodesic is strictly less than $\pi$ is dense in the space of all measured geodesic laminations.

\sk{Conclusion}\label{conclu}
It remains to glue altogether the main propositions of the previous sections to obtain the proof of our main result.

From the stability property, it suffices to build a transverse invariant measure of full support on ${\cal N}$, whenever $N$ has constant curvature. Let $\G=\pi_{1}(N)$, and $F=PSL(2,\mathbb Z)$.

We consider the restricted configuration space $\bbi$, subset of the space of maps from
$\mathbb {QP}^{1}$ to $\cpi$, as defined in paragraph \ref{configspa}, and  associated to the Markov configuration data (as defined in \ref{data}) coming from  the complex crossratio on $\di\yp$, according to section \ref{crossratio} and proposition
\ref{virdata}. 

We can turn now this configuration data into a measured one as shown in paragraph \ref{mesconfdata}. 

Thanks now to the main result of section \ref{combmod}, theorem \ref{exismes}, we obtain a finite $F$-invariant ergodic measure of full support on $\bbi/\Gamma$. Here, the action of $F$ is by right composition.

Furthermore, choices of mutually singular equivariant families of measures
lead to mutually singular transversal measure.

Next, we  suspend the action of $F$ on  $\bbi$. Namely, 
namely we consider the Riemannian lamination.
$$
{\cal F}=(\bbi\times \hyp )/F,              ,
$$
where $F$ acts as an isometry group on $\hyp$ and diagonally on
$\bbi\times\hyp$.

The finite ergodic and $F$-invariant measure on $\bbi/\G$ gives rise to a
transversal $\G$-invariant and ergodic measure on ${\cal F}$ that we
call $\mu$.  

Finally proposition \ref{exisphi} defines a  map $\Phi$ from  $\cal F /\G$ to $\cal N$, which is leaf preserving, continuous with a dense image. Therefore, we can pushforward $\mu$ using $\Phi$ to obtain a transversal ergodic  finite measure  of full support.

\auteur

\vskip 1truecm
\begin{thebibliography}{99}
\bibitem{la}{F.~Labourie} 
\emph{Un lemme de Morse pour les surfaces convexes.} {A para{\^\i}tre {\`a} {\it
Invent. Math.}}
\bibitem{la2}{F.~Labourie}
\emph{Immersions isom{\'e}triques elliptiques et courbes pseudo-holomorphes.}~J. Diff. Geom. {\bf 30} (1989) pp 395-44. 
\bibitem{la3}{F.~Labourie}
\emph{Surfaces convexes dans
l'espace hyperbolique et $\cpi$-structures.}~J. Lond. Mat. Soc. {\bf 45} (1992)pp 549-565.
\bibitem{gafa}{F.~Labourie}
\emph{Probl{\`e}mes de
Monge-Amp{\`e}res, courbes pseudo-holomorphes et laminations}~G.A.F.A. {\bf 7} (1997), pp 496-534.
\bibitem{led}{F.~Ledrappier}
\emph{Structure au bord {\`a} des vari{\'e}t{\'e}s {\`a} courbure
n{\'e}\-ga\-tive.}
S{\'e}\-mi\-nai\-re de th{\'e}orie spectrale et g{\'e}ometrie (1995) pp 97-122
\bibitem{Su}{D.~Sullivan}
\emph{The Density at Infinity of a Discrete Group of Hyperbolic Motions.}{~Public. Math. I.H.E.S. {\bf 50} (1979) pp 171-202}
\bibitem{tan}{H.~Tanigawa}
\emph{Grafting, Harmonic Maps and  Projective Structures on Surfaces.}{~J. of Diff. Geom. {\bf 47} (1997) pp 399-419}
\end{thebibliography}
\end{document}